\documentclass[preprint,review,12pt]{elsarticle}
\usepackage{amssymb}
\usepackage{amsfonts}
\usepackage{amsmath}
\usepackage{amsthm}
\usepackage{color}
\usepackage{graphicx}
\usepackage{epsfig,mathrsfs}
\usepackage[bf,SL,BF]{subfigure}
\usepackage{fancyhdr}
\usepackage{CJK}
\usepackage{caption}
\usepackage{wrapfig}
\usepackage{cases}
\usepackage{subfigure}
\graphicspath{{figure/}}

\newtheorem*{lemma*}{Lemma A}

\numberwithin{equation}{section}

\renewcommand{\eqref}[1]{(\ref{#1})}

 \allowdisplaybreaks

\begin{document}
 \pagenumbering{arabic}

\begin{frontmatter}
\title{{\bf \noindent A second-order numerical method for two-dimensional two-sided space fractional convection diffusion equation}}
\author{Minghua Chen, Weihua Deng$^{*}$}
\cortext[cor2]{Corresponding author. E-mail: dengwh@lzu.edu.cn.}
\address{School of Mathematics and Statistics,
Lanzhou University, Lanzhou 730000, P. R. China}

\date{}

\begin{abstract}

Space fractional convection diffusion equation describes physical
phenomena where particles or energy (or other physical quantities)
are transferred inside a physical system due to two processes:
convection and superdiffusion. In this paper, we discuss the
practical alternating directions implicit method to solve the
two-dimensional two-sided space fractional convection diffusion
equation  on a finite domain. We
theoretically prove and numerically verify that the presented finite
difference scheme is unconditionally von Neumann stable and second order
convergent in both space and time directions.

\end{abstract}
\end{frontmatter}

\section{Introduction}
Relaxing the restriction of the boundedness of the second moments,
fractional derivatives naturally appear to characterize anomalous
diffusion \cite{Metzler:00}, usually the time fractional derivative
is used for describing the subdiffusion and space fractional
derivative for the superdiffusion. More often, diffusion corresponds
to a power law, $\langle x^2(t) \rangle \sim Dt^\alpha$, where $D$
is the diffusion coefficient and $t$ the elapsed time. In a
classical diffusion process, $\alpha=1$. If $\alpha<1$, the
phenomenon is called subdiffusion. If $\alpha>1$, the particles
undergo superdiffusion. For $\alpha=1$, the position probability
density of the particles satisfies the classical diffusion equation
with first order time derivative and second order space derivative;
for $\alpha<1$, the position probability density of the particles is
the solution of the time fractional diffusion equation with
$\alpha$-th order time derivative and second order space derivative;
for $\alpha>1$, the position probability density satisfies the space
fractional diffusion equation with first order time derivative and
$(3-\alpha)$-th order space derivative.

The space fractional advection diffusion equation describes the
physical phenomena involving two physical processes: convection and
superdiffusion, i.e., for the equation, besides the $\alpha$-th
order space fractional derivative term there exists the classical
first order space derivative. This paper focuses on the numerical
algorithm of the following two-dimensional two-sided space
fractional convection diffusion equation
\begin{equation}\label{1.1}
\begin{split}
  \frac{\partial u(x,y,t) }{\partial t}=&d_+\left(x,y\right)~ _{x_L}D_x^{\alpha}u(x,y,t)+d_-(x,y)~  _{x}D_{x_R}^{\alpha}u(x,y,t) \\
  &+e_+(x,y) ~ _{y_L}D_y^{\beta}u(x,y,t)+e_-(x,y) ~  _{y}D_{y_R}^{\beta}u(x,y,t)\\
  &+g(x,y) u_x(x,y,t)+h(x,y) u_y(x,y,t),
\end{split}
\end{equation}
where $(x,y) \in \Omega=(x_L,x_R) \times (y_L,y_R),\, 0< t \leq
T$, the fractional orders $1<\alpha,\beta<2$; the function $
s(x,y,t)$ is a source term; and the diffusion coefficients
$d_+(x,y)\geq 0, d_-(x,y)\geq0, e_+(x,y)\geq0$ and $ e_-(x,y)\geq0$.
The initial and boundary conditions are, respectively, taken as
\begin{equation}
u(x,y,0)=u_0(x,y) ~~~ {\rm for}~~~ (x,y) \in \Omega,
\end{equation}
and
\begin{equation}
  u(x,y,t)|_{\partial \Omega} =B(x,y,t).
\end{equation}
The left  and  right  Riemann-Liouville fractional
derivatives of order $\mu$  ($0\leq n-1\leq \mu < n$ and $n$ is an integer)
are, respectively, defined by \cite{Kenneth:93,Podlubny:99,Samko:93}
\begin{equation}\label{1.4}
  _{x_L}D_x^{\mu}u(x)=D^n{[_{x_L}D_{x}^{-(n-\mu)}u(x)]},
\end{equation}
and
\begin{equation}\label{1.5}
_{x}D^{\mu}_{x_R}u(x)=E^n{[_{x}D_{x_R}^{-(n-\mu)}u(x)]},
\end{equation}
where
\begin{equation*}
   _{x_L}D^{-\nu}_{x}u(x)=\frac{1}{\Gamma (\nu)}\int_{x_L}^{x}{(x-\xi)^{\nu-1}}{u(\xi)}d\xi, \quad  \nu>0,
\end{equation*}

\begin{equation*}
   _{x}D^{-\nu}_{x_R}u(x)=\frac{1}{\Gamma (\nu)}\int_x^{x_R}{(\xi-x)^{\nu-1}}{u(\xi)}d\xi, \quad  \nu>0,
\end{equation*}
and
\begin{equation*}
   D=\frac{d}{dx},  \quad    \quad   E\equiv -D=-\frac{d}{dx}.
\end{equation*}
 The
Gr\"{u}nwald-Letnikov definitions for the left and right
 fractional derivatives are, respectively, given as
\begin{equation}\label{1.41}
{_{x_L}^{G}D_{x}^{\mu}}u(x)=\lim_{M_+\rightarrow
\infty}\frac{1}{h_+^\mu}{}\sum_{i=0}^{M_+}(-1)^i {\mu \choose
i}u(x-ih_+),
\end{equation}
and
\begin{equation}\label{1.51}
{_{x}^{G}D}_{x_R}^{\mu}u(x)=\lim_{M_-\rightarrow
\infty}\frac{1}{h_-^\mu}{}\sum_{i=0}^{M_-}(-1)^i {\mu \choose
i}u(x+ih_-),
\end{equation}
where $h_+=(x-x_L)/M_+,\,h_-=(x_R-x)/M_-$, and $M_+$ and $M_-$ are
positive integers. The Riemann-Liouville and Gr\"{u}nwald-Letnikov
derivatives are equivalent  \cite{Li:07} under the assumptions that
the function performed are sufficiently smooth.

Using the formulae $(\ref{1.41})$ and $(\ref{1.51})$ to discretize
the space fractional derivatives is a nature idea for designing the
numerical schemes. Unfortunately, this usually leads to
unconditionally unstable finite difference schemes, but Meerschaert
et al successfully circumvent this difficulties by modifying the
formulae $(\ref{1.41})$ and $(\ref{1.51})$ to obtain the so-called
shifted Gr\"{u}wald formulae \cite{Meerschaert:04}. Based on these
shifted formulae, Meerschaert et al did a series of works for
numerically solving space fractional diffusion equations
\cite{Meerschaert:06,Mark:06,Charles:07,Charles:06}. Sousa proposes
another way to approximate the fractional Caputo derivatives
\cite{Sousa:Pro,Sousa:11}, which can obtain second order accuracy.
More recently Sousa et al further discuss using the similar idea to
discretize Riemann-Liouville fractional derivatives in infinite
domain \cite{Sousa:12}. By using different ideas, Liu's group discusses the finite difference methods for fractional partial differential equations with Riesz space fractional derivatives \cite{Yang:10} and the
difference methods for the space-time fractional advection-diffusion equation \cite{Liu:07}.  For the discretization of the time
fractional derivatives, usually the different challenges will be met
\cite{Deng:07,Diethelm:02,Li:11,Liu:03}.
   Here we will combine the
alternating directions implicit (ADI) method with Crank-Nicolson
scheme to design the finite difference scheme for the two-sided
two-dimensional space fractional advection diffusion equation
(\ref{1.1}).  The numerical scheme will be theoretically proven and
numerically verified to be unconditionally von Neumann stable and
second order convergent.

The  paper is organized  as follows. In Section 2, we derive the
linear spline approximation to the right Riemann-Liouville
fractional derivative, and the full discretization of (\ref{1.1}) is
presented, where the Crank-Nicolson scheme and the alternating
directions implicit  method are combined together. Section 3 does
the detailed theoretical analyses for the consistency and stability
of the given schemes. To show the effectiveness of the algorithm, we
perform the numerical experiments to verify the theoretical results
 in Section 4. Finally, we
conclude the paper with some remarks in the last section.

\section{Discretization Schemes}\label{sec:1}

We use three subsections to derive the full discretization of
(\ref{1.1}). Since the linear spline approximation for the left
Riemann-Liouville fractional derivative in the infinite interval can
be easily got from \cite{Sousa:12}, we further derive the linear
spline approximation for the right Riemann-Liouville fractional
derivative and make some remarks on the relationship between the
discretization schemes of left and right Riemann-Liouville
fractional derivatives in the first subsection. Then in the second
subsection, we present the scheme for the one dimensional case of
(\ref{1.1}). The third subsection detailedly provides the full
discrete scheme of the two-dimensional two-sided space fractional
convection diffusion equation.

\subsection{Discretizations for the left and right Riemann-Liouville fractional derivatives}
Let the mesh points $x_i=x_L+i\Delta x,i=0,1,\ldots ,{N_x},$ where
$\Delta x=(x_R-x_L)/{N_x}$ is the uniform space step. Taking
$\alpha\in (1,2)$ in the left Riemann-Liouville fractional
derivative (\ref{1.4}), its approximation operator
$\delta'_{\alpha,_+x}{u_i^n}$ has second order accuracy in a bounded
domain (proved in next section), where ${u_i^n}$ denotes the
approximated value of $u(x_i,t_n)$ and the left fractional
approximation operator, which can be obtained by truncating its
infinite version \cite{Sousa:12}, is defined as
\begin{equation}\label{2.19}
  \delta'_{\alpha,_+x}{u_i^n}:=\frac{1}{\Gamma(4-\alpha)\Delta
  x^{\alpha}}\sum_{k=0}^{i+1}u_k^np_{i,k}^{\alpha},
\end{equation}
where
\begin{equation}\label{2.20}
p_{i,k}^{\alpha}=\left\{ \begin{array}
 {l@{\quad } l}
  a_{i-1,k}-2a_{i,k}+a_{i+1,k},&k \leq i-1,\\
  -2a_{i,i}+a_{i+1,i} ,&k=i,\\
 a_{i+1,i+1} , & k=i+1,\\
 0 , &k>i+1,
 \end{array}
 \right.
\end{equation}
and
\begin{equation*}
 a_{i,k}=\left\{ \begin{array}
 {l@{\quad} l}
 (i-1)^{3-\alpha}-i^{2-\alpha}(i-3+\alpha),& k=0,\\
  (i-k+1)^{3-\alpha}-2(i-k)^{3-\alpha}+(i-k-1)^{3-\alpha},&1 \leq k \leq i-1,\\
  1,&k=i.
 \end{array}
 \right.
\end{equation*}

Here in the finite interval $x_L< x <x_R$ with $\alpha \in (1,2)$,
we further derive the linear spline approximation for the right
Riemann-Liouville fractional derivative defined by
%
\begin{equation} \label{1.12}
   _{x}D_{x_R}^{\alpha}u(x,t)=\frac{1}{\Gamma (2-\alpha)}\frac{\partial^2}{\partial x^2}
  \int_x^{x_R} u(\xi,t)(\xi-x)^{1-\alpha}d\xi.
\end{equation}
For a fixed time t, denote
\begin{equation}
\mathcal{ I }_{\alpha}(x)=\frac{1}{\Gamma (2-\alpha)}\int_x^{x_R} u(\xi,t)(\xi-x)^{1-\alpha}d\xi,
\end{equation}
then
\begin{equation}\label{1.8}
  _{x}D_{x_R}^{\alpha}u(x,t)=\frac{\partial^2}{\partial x^2}
  \mathcal{ I }_{\alpha}(x),
\end{equation}
and we can do the following approximation at $x_i$,
\begin{equation}\label{2.10}
  \frac{\partial ^2}{\partial x^2}\mathcal{ I }_{\alpha}(x_i)\simeq\frac{1}{\Delta x^2}
  [\mathcal{ I }_{\alpha}(x_{i-1})-2\mathcal{ I }_{\alpha}(x_{i})+\mathcal{ I
  }_{\alpha}(x_{i+1})],~~~1 \leq i \leq N_x-1.
\end{equation}
For each $x_i$, take
\begin{equation}\label{1.11}
 \mathcal{I}_{\alpha}(x_i)\simeq I_{\alpha}(x_i)=\frac{1}{\Gamma (2-\alpha)}\int_{x_i}^{x_R}
 S_i(\xi)(\xi-x_i)^{1-\alpha}d\xi,
\end{equation}
where  the spline $S_i(\xi)$ is defined by
\begin{equation}\label{2.12}
  S_i(\xi)=\sum_{k=i}^{N_x}u(x_k,t)s_{i,k}(\xi),
\end{equation}
with $s_{i,k}(\xi)$, in every subinterval $[x_{k-1},x_{k+1}]$, for
$i+1\leq k \leq {N_x}-1$, given as
\begin{equation*}
 s_{i,k}(\xi)=\left\{ \begin{array}
 {l@{\quad } l}
 \displaystyle\frac{\xi-x_{k-1}}{x_{k}-x_{k-1}},  &x_{k-1}\leq \xi \leq x_{k}, \\
 \\
   \displaystyle\frac{x_{k+1}-\xi}{x_{k+1}-x_{k}}, &x_{k}\leq \xi \leq x_{k+1}, \\
   \\
  0,& {\rm otherwise},
 \end{array}
 \right.
\end{equation*}
and for $k=i$ and $k={N_x}$, $s_{i,k}(\xi)$ taken as
\begin{equation*}
 s_{i,i}(\xi)=\left\{ \begin{array}
 {l@{\quad} l}
  \displaystyle\frac{x_{i+1}-\xi}{x_{i+1}-x_{i}},&x_{i}\leq \xi \leq x_{i+1}, \\
  \\
  0,& {\rm otherwise},
 \end{array}
 \right.
\end{equation*}
and
\begin{equation*}
 s_{i,{N_x}}(\xi)=\left\{ \begin{array}
 {l@{\quad } l}
 \displaystyle \frac{\xi-x_{{N_x}-1}}{x_{N_x}-x_{N_x-1}},&x_{N_x-1}\leq \xi \leq x_{{N_x}}, \\
 \\
 0,& {\rm otherwise}.
 \end{array}
 \right.
\end{equation*}
 According to (\ref{1.11}) and (\ref{2.12}), we have
\begin{equation}
\begin{split}
 \mathcal{I}_{\alpha}(x_i)\simeq I_{\alpha}(x_i)
  &=\frac{1}{\Gamma (2-\alpha)}u(x_i,t)\int_{x_{i}}^{x_{i+1}} s_{i,i}(\xi)(\xi-x_i)^{1-\alpha}d\xi\\
 &\quad +\frac{1}{\Gamma (2-\alpha)}\sum_{k=i+1}^{N_x-1}u(x_k,t)\int_{x_{k-1}}^{x_{k+1}} s_{i,k}(\xi)(\xi-x_i)^{1-\alpha}d\xi\\
 &\quad +\frac{1}{\Gamma (2-\alpha)}u(x_{N_x},t)\int_{x_{{N_x}-1}}^{x_{N_x}} s_{i,{N_x}}(\xi)(\xi-x_i)^{1-\alpha}d\xi,
 \end{split}
\end{equation}
and
\begin{equation}\label{2.10905}
\begin{split}
 \frac{1}{\Gamma (2-\alpha)}\int_{x_{i}}^{x_{i+1}} s_{i,i}(\xi)(\xi-x_i)^{1-\alpha}d\xi &=\frac{\Delta x^{2-\alpha}}{\Gamma(4-\alpha)}b_{i,i},\\
\frac{1}{\Gamma (2-\alpha)}\int_{x_{k-1}}^{x_{k+1}} s_{i,k}(\xi)(\xi-x_i)^{1-\alpha}d\xi& =\frac{\Delta x^{2-\alpha}}{\Gamma(4-\alpha)}b_{i,k},\\
\frac{1}{\Gamma (2-\alpha)}\int_{x_{{N_x}-1}}^{x_{N_x}} s_{i,{N_x}}(\xi)(\xi-x_i)^{1-\alpha}d\xi& =\frac{\Delta x^{2-\alpha}}{\Gamma(4-\alpha)}b_{i,N_x},\\
 \end{split}
\end{equation}
where
\begin{equation}\label{2.105}
 b_{i,k}=\left\{ \begin{array}
 {l@{\quad } l}
   1,&k=i, \\
  (k-i+1)^{3-\alpha}-2(k-i)^{3-\alpha}+(k-i-1)^{3-\alpha},&i+1 \leq k \leq {N_x}-1,\\
 (3-\alpha-{N_x}+i)({N_x}-i)^{2-\alpha}+({N_x}-i-1)^{3-\alpha},& k={N_x}.
 \end{array}
 \right.
\end{equation}
Then
\begin{equation}\label{2.155}
 \mathcal{I}_{\alpha}(x_i)\simeq I_{\alpha}(x_i)=\frac{\Delta x^{2-\alpha}}{\Gamma(4-\alpha)}\sum_{k=i}^{N_x}u(x_k,t)b_{i,k},
\end{equation}
and (\ref{2.10}) can be written as
\begin{equation}\label{2.16}
\begin{split}
  \frac{\partial ^2}{\partial x^2}\mathcal{ I }_{\alpha}(x_i) &\simeq \frac{1}{\Delta x^2}
  [\mathcal{ I }_{\alpha}(x_{i-1})-2\mathcal{ I }_{\alpha}(x_{i})+\mathcal{ I }_{\alpha}(x_{i+1})]\\
  &\simeq\frac{1}{\Gamma(4-\alpha)\Delta x^{\alpha}}\left[\sum_{k=i-1}^{N_x}u(x_k,t)b_{i-1,k}-2\sum_{k=i}^{N_x}u(x_k,t)b_{i,k}+\sum_{k=i+1}^{N_x}u(x_k,t)b_{i+1,k}\right ]\\
  &=\frac{1}{\Gamma(4-\alpha)\Delta x^{\alpha}}\sum_{k=i-1}^{N_x}u(x_k,t)q_{i,k}^{\alpha},
  \end{split}
\end{equation}
where
\begin{equation}\label{2.17}
q_{i,k}^{\alpha}=\left\{ \begin{array}
 {l@{\quad } l}
   0,&k<i-1, \\
  b_{i-1,i-1},&k=i-1,\\
 -2b_{i,i}+b_{i-1,i}, & k=i,\\
 b_{i-1,k}-2b_{i,k}+b_{i+1,k}, &i+1\leq k \leq {N_x}.
 \end{array}
 \right.
\end{equation}
Denoting ${u_i^n}$ as the approximated value of $u(x_i,t_n)$, we can
define the right fractional approximation operator as
\begin{equation}\label{2.18}
  \delta'_{\alpha,_-x}{u_i^n}:=\frac{1}{\Gamma(4-\alpha)\Delta
  x^{\alpha}}\sum_{k=i-1}^{N_x}u_k^nq_{i,k}^{\alpha},
\end{equation}
which has second order accuracy for approximating (\ref{1.12})
(proved in the next section).

\noindent{\bf Remark 2.1.} Denoting $\tilde{U}^n=[u_1^n,u_2^n,\cdots,u_{N_x-1}^n]^{\rm T}$, and
 rewriting (\ref{2.19}) and (\ref{2.18}) as matrix forms
$\delta'_{\alpha,_+x} \tilde{U}^n=\tilde{A}\tilde{U}^n+b_1$ and
$\delta'_{\alpha,_-x} \tilde{U}^n=\tilde{B}\tilde{U}^n+b_2$,
respectively, then there exists $\tilde{A}=\tilde{B}^{\rm T}$.

\subsection{Numerical scheme for one-dimensional fractional convection diffusion equation}
   We now examine the full discretization scheme to the one-dimensional two-sided fractional convection diffusion equation
\begin{equation}\label{2.21}
\begin{split}
  \frac{\partial u(x,t) }{\partial t}=d_+(x)~ _{x_L}D_{x}^{\alpha}u(x)+d_-(x)~ _{x}D_{x_R}^{\alpha}u(x)
  +g(x) u_x(x,t)+s(x,t).
\end{split}
\end{equation}
In the time direction, we use the Crank-Nicolson scheme.  The
central difference formula, left fractional approximation operator
(\ref{2.19}), and right fractional approximation operator
(\ref{2.18}) are respectively used to discretize the classical first
order space derivative, left Riemann-Liouville fractional
derivative, and right Riemann-Liouville fractional derivative.
Taking the uniform time step $\Delta t$ and space step $\Delta x$,
and setting $ d_{+,i}=d_+(x_i), d_{-,i}=d_-(x_i), g_i=g(x_i), $ and
$s_i^{n+1/2}=s(x_i,t_{n+1/2})$, where $t_{n+1/2}=(t_n+t_{n+1})/2$,
the full discretization of (\ref{2.21}) has the following form
\begin{equation}\label{2.22}
\begin{split}
  \frac{u_i^{n+1}-u_i^n}{\Delta t}=&\frac{1}{\Gamma(4-\alpha)\Delta x^{\alpha}}\left [ \sum_{k=0}^{i+1}p_{i,k}^{\alpha}d_{+,i}\frac{u_k^{n+1}+u_k^{n}}{2}+\sum_{k=i-1}^{N_x}q_{i,k}^{\alpha}d_{-,i}\frac{u_k^{n+1}+u_k^{n}}{2}\right ] \\
  &+\frac{g_i}{2\Delta x} \left ( \frac{u_{i+1}^{n+1}+u_{i+1}^n}{2}- \frac{u_{i-1}^{n+1}+u_{i-1}^n}{2} \right )
  +s_i^{n+1/2}.
\end{split}
\end{equation}
Similar to (\ref{2.19}) and (\ref{2.18}), we define
\begin{equation} \label{2.23}
\begin{split}
& D'_{\alpha,x}u_i^n:=\frac{u_{i+1}^{n}-u_{i-1}^n}{2\Delta x};  \\
&D''_{\alpha,x}u_i^n:=\frac{u_{i+1}^{n}-u_{i-1}^n}{2\Delta x}g_i;\\
& { \delta}_{\alpha,_+x}^{''}{u_i^n}:=\frac{d_{+,i}}{\Gamma(4-\alpha)\Delta x^{\alpha}}\sum_{k=0}^{i+1}u_k^np_{i,k}^{\alpha};\\
&{
\delta}_{\alpha,_-x}^{''}{u_i^n}:=\frac{d_{-,i}}{\Gamma(4-\alpha)\Delta
x^{\alpha}}\sum_{k=i-1}^{{N_x}}u_k^nq_{i,k}^{\alpha},
 \end{split}
 \end{equation}
 then  (\ref{2.22}) can be expressed as
\begin{equation}\label{2.26}
 \left [1-\frac{\Delta t}{2}\left( \delta''_{\alpha,_+x} +  \delta''_{\alpha,_-x} + D''_{\alpha,x}   \right) \right ]u_i^{n+1}=
 \left [1+\frac{\Delta t}{2} \left( \delta''_{\alpha,_+x} +  \delta''_{\alpha,_-x} + D''_{\alpha,x}  \right ) \right ]u_i^{n}+  s_i^{n+1/2}\Delta t,
\end{equation}
for $i=1,2,\ldots,{N_x}-1$, associated with the boundary conditions
$u_0^n$ and $u_{N_x}^n$.  Putting $\xi_i=\frac{\Delta
t}{2\Gamma\left(4-\alpha\right)\Delta x^{\alpha}}d_{+,i}$,
$\eta_i=\frac{\Delta t}{2\Gamma(4-\alpha)\Delta x^{\alpha}}d_{-,i},$
and $\gamma_i=\frac{\Delta t}{4\Delta x}g_i$, the system of
equations given by (\ref{2.22}) takes the form
   \begin{equation}(I-A)U^{n+1}=(I+A)U^{n}+\Delta t S^{n+1/2},\end{equation}
 where $I$ is the identity  matrix, and
 \begin{equation*}
 U^{n}=[u_0^n,u_1^n,u_2^n,\ldots,u_{N_x}^n]^{\rm T}, ~~S^{n+1/2}=[s_0^{n+1/2},s_1^{n+1/2},s_2^{n+1/2},\ldots,s_{N_x}^{n+1/2}]^{\rm T},
  \end{equation*}
 and the matrix entries $A_{i,j}$ for $i=1,\ldots,N_x-1$ and $j=1,\ldots,N_x-1 $
  are defined by
  \begin{equation}\label{c28}
A_{i,j}=\left\{ \begin{array}
 {l@{\quad} l}
 \xi_ip_{i,i}^{\alpha}+\eta_iq_{i,i}^{\alpha},& j=i,\\
\xi_ip_{i,i-1}^{\alpha}+\eta_iq_{i,i-1}^{\alpha}+\gamma_i,&j=i-1,\\
\xi_ip_{i,i+1}^{\alpha}+\eta_iq_{i,i+1}^{\alpha}-\gamma_i,&j=i+1,\\
\xi_ip_{i,j}^{\alpha},&j<i+1,\\
\eta_iq_{i,j}^{\alpha},&j>i+1,\\
 \end{array}
 \right.
\end{equation}
and $A_{0,0}=1,A_{0,i}=0$ for $i=1,\ldots,N_x;\,A_{N_x,N_x}=1$ and
$A_{N_x,i}=0$ for $ i=0,\ldots, N_x-1$.

%
%
\subsection{ADI scheme for two-dimensional two-sided fractional convection diffusion equation}
Under the direction of discretizing the one-dimensional case of
(\ref{1.1}) in the last subsection, we use ADI
\cite{Charles:07,Gustafsson:95,Lapidus:82} to numerically solve
(\ref{1.1}). First we introduce and list the denotations (some of
them already given above) that will be used in the following:
\begin{equation}\label{2.28}
\begin{split}
&D'_{\alpha,x}u_{i,j}^n:=\frac{u_{i+1,j}^{n}-u_{i-1,j}^n}{2\Delta
x};
\qquad\qquad\qquad\qquad \, D'_{\beta,y}u_{i,j}^n:=\frac{u_{i,j+1}^{n}-u_{i,j-1}^n}{2\Delta y};\\
&D''_{\alpha,x}u_{i,j}^n:=\frac{u_{i+1,j}^{n}-u_{i-1,j}^n}{2\Delta
x}g_{i,j};
\qquad\qquad \qquad\quad  D''_{\beta,y}u_{i,j}^n:=\frac{u_{i,j+1}^{n}-u_{i,j-1}^n}{2\Delta y}h_{i,j};\\
&\delta'_{\alpha,_+x}{u_{i,j}^n}:=\frac{1}{\Gamma(4-\alpha)\Delta
x^{\alpha}}\sum_{k=0}^{i+1}u_{k,j}^np_{i,k}^{\alpha};
\qquad\quad \,\,\,\, \delta'_{\beta,_+y}{u_{i,j}^n}:=\frac{1}{\Gamma(4-\beta)\Delta y^{\alpha}}\sum_{k=0}^{j+1}u_{i,k}^np_{j,k}^{\beta};\\
& \delta'_{\alpha,_-x}{u_{i,j}^n}:=\frac{1}{\Gamma(4-\alpha)\Delta
x^{\alpha}}\sum_{k=i-1}^{N_x}u_{k,j}^nq_{i,k}^{\alpha};
\qquad \,\,\,\,\,\,\delta'_{\beta,_-y}{u_{i,j}^n}:=\frac{1}{\Gamma(4-\beta)\Delta y^{\alpha}}\sum_{k=j-1}^{N_y}u_{i,k}^nq_{j,k}^{\beta};\\
&{
\delta''}_{\alpha,_+x}{u_{i,j}^n}:=\frac{d_{+,i,j}}{\Gamma(4-\alpha)\Delta
x^{\alpha}}\sum_{k=0}^{i+1}u_{k,j}^np_{i,k}^{\alpha};
\qquad\quad\,\,\,\, \,{ \delta''}_{\beta,_+y}{u_{i,j}^n}:=\frac{e_{+,i,j}}{\Gamma(4-\beta)\Delta y^{\alpha}}\sum_{k=0}^{j+1}u_{i,k}^np_{j,k}^{\beta};\\
&{
\delta''}_{\alpha,_-x}{u_{i,j}^n}:=\frac{d_{-,i,j}}{\Gamma(4-\alpha)\Delta
x^{\alpha}}\sum_{k=i-1}^{{N_x}}u_{k,j}^nq_{i,k}^{\alpha};
\qquad\,\,\,\,\,\,\,{ \delta''}_{\beta,_-y}{u_{i,j}^n}:=\frac{e_{-,i,j}}{\Gamma(4-\beta)\Delta y^{\alpha}}\sum_{k=j-1}^{{N_y}}u_{i,k}^nq_{j,k}^{\beta}.\\
\end{split}
\end{equation}
 Analogously we still use the Crank-Nicolson scheme to do the discretization in time direction.
 Taking $u_{i,j}^n$ as the approximated value of $u(x_i,y_j,t_n)$,
 $d_{+,i,j}=d_+(x_i,y_j),\,d_{-,i,j}=d_-(x_i,y_j),\,e_{+,i,j}=e_+(x_i,y_j)$,
 $e_{-,i,j}=e_-(x_i,y_j),\,g_{i,j}=g(x_i,y_j),\,h_{i,j}=h(x_i,y_j)$, $ t_{n+1/2}=(t_n+t_{n+1})/2$, $s_{i,j}^{n+1/2}=s(x_i,y_j,t_{n+1/2})$,
 $\Delta x=(x_R-x_L)/{N_x}$, and
$\Delta y=(y_R-y_L)/{N_y}$, for the uniform space steps $\Delta
x,\Delta y$ and time  step $\Delta t$, the resulting discretization
of (\ref{1.1}) can be written as
\begin{equation}
\begin{split}
  \frac{u_{i,j}^{n+1}-u_{i,j}^n}{\Delta t}
  =&\frac{1}{\Gamma(4-\alpha)\Delta x^{\alpha}} \left[\sum_{k=0}^{i+1}p_{i,k}^{\alpha}d_{+,i,j}\frac{u_{k,j}^{n+1}+u_{k,j}^{n}}{2}+
     \sum_{k=i-1}^{N_x}q_{i,k}^{\alpha}d_{-,i,j}\frac{u_{k,j}^{n+1}+u_{k,j}^{n}}{2}\right ] \\
  &+\frac{1}{\Gamma(4-\beta)\Delta y^{\beta}} \left[\sum_{k=0}^{j+1}p_{j,k}^{\beta}e_{+,i,j}\frac{u_{i,k}^{n+1}+u_{i,k}^{n}}{2}+
     \sum_{k=j-1}^{N_y}q_{j,k}^{\beta}e_{-,i,j}\frac{u_{i,k}^{n+1}+u_{i,k}^{n}}{2}\right ] \\
  &+\frac{g_{i,j}}{2\Delta x} \left ( \frac{u_{i+1,j}^{n+1}+u_{i+1,j}^n}{2}- \frac{u_{i-1,j}^{n+1}+u_{i-1,j}^n}{2} \right )\\
  &+\frac{h_{i,j}}{2\Delta y} \left ( \frac{u_{i,j+1}^{n+1}+u_{i,j+1}^n}{2}- \frac{u_{i,j-1}^{n+1}+u_{i,j-1}^n}{2} \right )\\
  &+s_{i,j}^{n+1/2}.
\end{split}
\end{equation}
Using the notations (\ref{2.28}), we have
\begin{equation}\label{2.30}
\begin{split}
 &\left [1-\frac{\Delta t}{2}\left( \delta''_{\alpha,_+x}+  \delta''_{\alpha,_-x} + D''_{\alpha,x}  \right )
 -\frac{\Delta t}{2}\left( \delta''_{\beta,_+y} +  \delta_{\beta,_-y}^{''} + D''_{\beta,y}   \right)\right ]u_{i,j}^{n+1} \\
 &\quad=\left [1+\frac{\Delta t}{2}\left( \delta''_{\alpha,_+x} +  \delta''_{\alpha,_-x} + D''_{\alpha,x}  \right )
 +\frac{\Delta t}{2}\left( \delta''_{\beta,_+y} +  \delta''_{\beta,_-y} + D''_{\beta,y} \right )   \right ]u_{i,j}^{n}+  s_{i,j}^{n+1/2}\Delta t.
\end{split}
\end{equation}
Further define
\begin{equation}
\begin{split}
~~~\delta_{\alpha,x}:= \delta''_{\alpha,_+x} +  \delta''_{\alpha,_-x} + D''_{\alpha,x};   \\
 \delta_{\beta,y}:=\delta''_{\beta,_+y} +  \delta''_{\beta,_-y} +
 D''_{\beta,y},
\end{split}
\end{equation}
thus, (\ref{2.30}) may be rewritten as
\begin{equation}\label{2.32}
\begin{split}
 &\left (1-\frac{\Delta t}{2}\delta_{\alpha,x}-\frac{\Delta t}{2}\delta_{\beta,y}\right )u_{i,j}^{n+1} \\
 &=\left (1+\frac{\Delta t}{2}\delta_{\alpha,x}+\frac{\Delta t}{2}\delta_{\beta,y}\right )u_{i,j}^{n}+  s_{i,j}^{n+1/2}\Delta t. \\
\end{split}
\end{equation}
For the two-dimensional two-sided fractional convection diffusion
equation (\ref{1.1}), the relevant perturbation of (\ref{2.32}) is
of the form
\begin{equation}\label{2.33}
\begin{split}
 &\left (1-\frac{\Delta t}{2}\delta_{\alpha,x} \right )\left(1-\frac{\Delta t}{2}\delta_{\beta,y}\right )u_{i,j}^{n+1} \\
 &=\left (1+\frac{\Delta t}{2}\delta_{\alpha,x}\right)\left(1+\frac{\Delta t}{2}\delta_{\beta,y}\right )u_{i,j}^{n}+  s_{i,j}^{n+1/2}\Delta t. \\
\end{split}
\end{equation}
The system of equations defined by (\ref{2.33}) may be solved by the
following  ADI  (Peaceman-Rachford type) scheme \cite{Charles:07}:
\begin{equation}\label{2.34}
\left (1-\frac{\Delta t}{2}\delta_{\alpha,x} \right )u_{i,j}^{*}
 =\left(1+\frac{\Delta t}{2}\delta_{\beta,y}\right )u_{i,j}^{n}+ \frac{\Delta t}{2}
 s_{i,j}^{n+1/2};
\end{equation}
\begin{equation}\label{2.35}
\left(1-\frac{\Delta t}{2}\delta_{\beta,y}\right )u_{i,j}^{n+1}
=\left (1+\frac{\Delta t}{2}\delta_{\alpha,x}\right)u_{i,j}^{*}+
\frac{\Delta t}{2} s_{i,j}^{n+1/2},
\end{equation}
where $u_{i,j}^{*}$ is an intermediate solution. For maintaining the
consistency, we need to carefully specify the boundary conditions of
$u_{i,j}^{*}$.
Subtracting (\ref{2.35}) from  ({\ref{2.34}), we obtain
\begin{equation}\label{2.36}
2u_{i,j}^{*}=\left (1-\frac{\Delta t}{2}\delta_{\beta,y} \right
)u_{i,j}^{n+1}
 +\left(1+\frac{\Delta t}{2}\delta_{\beta,y}\right )u_{i,j}^{n},
\end{equation}
then the boundary conditions  for $u_{i,j}^{*}$ ( $i=0$ and
$i={N_x}$ with $j=1,\ldots,{N_y}-1)$ can be given as
\begin{equation}\label{2.37}
  \begin{split}
    &u_{0,j}^{*}\,\,\,=\frac{1}{2}\left[\left (1-\frac{\Delta t}{2}\delta_{\beta,y} \right )u_{0,j}^{n+1}
 +\left(1+\frac{\Delta t}{2}\delta_{\beta,y}\right )u_{0,j}^{n}\right] \\
 &\qquad\,=\frac{1}{2}\left[\left (1-\frac{\Delta t}{2}\delta_{\beta,y} \right )B_{0,j}^{n+1}
 +\left(1+\frac{\Delta t}{2}\delta_{\beta,y}\right
 )B_{0,j}^{n}\right];
 \\
  &u_{{N_x},j}^{*}=\frac{1}{2}\left[\left (1-\frac{\Delta t}{2}\delta_{\beta,y} \right )u_{{N_x},j}^{n+1}
 +\left(1+\frac{\Delta t}{2}\delta_{\beta,y}\right )u_{{N_x},j}^{n}\right] \\
 &\qquad\,=\frac{1}{2}\left[\left (1-\frac{\Delta t}{2}\delta_{\beta,y} \right )B_{{N_x},j}^{n+1}
 +\left(1+\frac{\Delta t}{2}\delta_{\beta,y}\right )B_{{N_x},j}^{n}\right].
  \end{split}
\end{equation}

The corresponding  procedure is executed   as follows:
\begin{description}
 \item[(1)] First for every fixed  $y=y_k$ $ (k=1,\ldots,{N_y}-1)$, solving a set of ${N_x}-1$ equations defined by (\ref{2.34}) at the mesh points $ x_i,i=1,\ldots,{N_x}-1$, to get
 $u_{i,k}^{*}$;
 \item[(2)] Next alternating the spatial direction, and for each fixed  $x=x_k$ $(k=1,\ldots,{N_x}-1)$ solving a set of ${N_y}-1$  equations defined by (\ref{2.35}) at the points $y_j,j=1,\ldots,{N_y}-1$, to obtain $u_{k,j}^{n+1}$.
\end{description}

\section{Convergence and Stability Analysis}
We show the convergence for one-dimensional and two-dimensional
two-sided fractional convection diffusion equation by proving the
consistency and stability (according to Lax's equivalence theorem).
\subsection{Convergence and stability for one dimensional two-sided fractional convection diffusion equation}

\noindent{\bf Lemma 3.1.}
 Let $u\in C^4(\bar{\Omega} )$ and $\xi\in[x_k,x_{k+1}]$, then we have
 \begin{equation*}
 u(\xi)-S_i(\xi)= \displaystyle\theta(1-\theta)\frac{(\Delta x)^2}{2!}u''(\xi)+\frac{\theta(1-\theta)(1-2\theta)(\Delta x)^3}{3!}u'''(\xi) +\mathcal{O}\left((\Delta x)^4 \right),
 \end{equation*}
 where $\Delta x=x_{k+1}-x_k$, $\theta=(\xi-x_k)/\Delta x$.

\begin{proof}
For any $\xi\in[x_k,x_{k+1}]$, using Taylor series expansion at $\xi$, it is easy to get the results.

  \end{proof}

\noindent{\bf Lemma 3.2.} \cite{Kenneth:93,Podlubny:99}
Let $D^{n}u(x)$ be continuous in the interval $[x_L,x_R]$ and $0\leq n-1 \leq \mu <n$, then for $x_L < x < x_R$ the following holds
\begin{equation*}
  \begin{split}
_{x_L}\!D^{\mu}_xu(x)=D^n{[_{x_L}D_x^{-(n-\mu)}u(x)]}
=_{x_L}\!\!\!D^{-(n-\mu)}_x[D^n u(x)]+\sum _{k=0}^{n-1}\frac{(x-x_L)^{k-\mu}}{\Gamma{(-\mu+k+1)}}D^{k}f(x_L).
  \end{split}
\end{equation*}

\noindent{\bf Lemma 3.3.}
Let $D^{n}u(x)$ be continuous in the interval $[x_L,x_R]$ and $0\leq n-1 \leq \mu <n$, then for $x_L < x < x_R$ the following holds
\begin{equation*}
  \begin{split}
_{x}D^{\mu}_{x_R}u(x)=E^n{[_{x}D_{x_R}^{-(n-\mu)}u(x)]}
=_{x}\!\!D^{-(n-\mu)}_{x_R}[E^n u(x)]+\sum _{k=0}^{n-1}\frac{({x_R}-x)^{k-\mu}}{\Gamma{(-\mu+k+1)}}E^{k}f({x_R}).
  \end{split}
\end{equation*}

\begin{proof}

 Since $ _{x}D^{-\nu}_{x_R}u(x)=\frac{1}{\Gamma (\nu)}\int_x^{x_R}{(\xi-x)^{\nu-1}}{u(\xi)}d\xi$,  $\nu>0$,
 it can be obtained
 \begin{equation*}
   \begin{split}
 D[ _{x}D^{-\nu}_{x_R}u(x)]&=\frac{1}{\Gamma (\nu)}\int_x^{x_R}{\frac{d(\xi-x)^{\nu-1}}{dx}}{u(\xi)}d\xi
                           =-\frac{1}{\Gamma (\nu)}\int_x^{x_R}{u(\xi)}d(\xi-x)^{\nu-1}\\
                           &=  _{x}\!\!D^{-\nu}_{x_R}[Du(x)] - \frac{(x_R-x)^{\nu-1}}{\Gamma (\nu)}f(x_R),
   \end{split}
 \end{equation*}
then
   \begin{equation}\label{3.00001}
   \begin{split}
 E[ _{x}D^{-\nu}_{x_R}u(x)]=  _{x}\!\!D^{-\nu}_{x_R}[Eu(x)] + \frac{(x_R-x)^{\nu-1}}{\Gamma (\nu)}f(x_R).
   \end{split}
 \end{equation}
 Using $E$ to differentiate both sides of (\ref{3.00001}) leads to
     \begin{equation}\label{3.00002}
   \begin{split}
 E^2[ _{x}D^{-\nu}_{x_R}u(x)]=E\{E[ _{x}D^{-\nu}_{x_R}u(x)]\}= E\{ _{x}D^{-\nu}_{x_R}[Eu(x)]\} + \frac{(x_R-x)^{\nu-2}}{\Gamma (\nu-1)}f(x_R).
   \end{split}
 \end{equation}
Replacing $u$ with $Eu$ in (\ref{3.00001}), then (\ref{3.00002}) can be rewritten as
\begin{equation}\label{3.00003}
   \begin{split}
 E^2[ _{x}D^{-\nu}_{x_R}u(x)]=  _{x}\!\!D^{-\nu}_{x_R}[E^2u(x)] + \frac{(x_R-x)^{\nu-1}}{\Gamma (\nu)}Ef(x_R) + \frac{(x_R-x)^{\nu-2}}{\Gamma (\nu-1)}f(x_R).
   \end{split}
 \end{equation}
 Repeated iterations establish the desired result.

\end{proof}

  \noindent{\bf Theorem 3.4.}\label{Theorem:5}
  Let  $u\in C^4(\bar{\Omega})$ and  $_xD_{x_R}^{\alpha-2}u\in C^4(\bar{\Omega})$,
   then $\displaystyle_{x}D_{x_R}^{\alpha}u(x_i,t_n)= \delta'_{\alpha,_-x}{u_i^n}+\mathcal{O}\left((\Delta
  x)^2\right)$, where $\alpha \in (1,2)$.

  \begin{proof}
 In the following we use $u(x_i)$ to denote
$u(x_i,t_n)$.  Since $t_n$ is fixed and  $_xD_{x_R}^{\alpha-2}u\in C^4(\bar{\Omega})$, then the following holds
 \begin{equation*}
 \begin{split}
  _{x}D_{x_R}^{\alpha}u(x_i)=\frac{\partial ^2}{\partial x^2}\mathcal{ I }_{\alpha}(x_i)=\frac{1}{(\Delta x)^2}
  [\mathcal{ I }_{\alpha}(x_{i-1})-2\mathcal{ I }_{\alpha}(x_{i})+\mathcal{ I }_{\alpha}(x_{i+1})]+\mathcal{O}\left((\Delta
  x)^2 \right).
  \end{split}
\end{equation*}
Let $\epsilon(x_i)$ be the error satisfying
\begin{equation*}
\mathcal{ I }_{\alpha}(x_{i-1})-2\mathcal{ I
}_{\alpha}(x_{i})+\mathcal{ I }_{\alpha}(x_{i+1})= I
_{\alpha}(x_{i-1})-2 I _{\alpha}(x_{i})+ I
_{\alpha}(x_{i+1})+\epsilon(x_i),
\end{equation*}
then there exists
\begin{equation*}
_{x}D_{x_R}^{\alpha}u(x_i)=\frac{1}{(\Delta x)^2}
  \left[I _{\alpha}(x_{i-1})-2I _{\alpha}(x_{i})+ I _{\alpha}(x_{i+1})\right] +\frac{1}{(\Delta x)^2}\epsilon(x_i)  +\mathcal{O}\left((\Delta
  x)^2\right),
\end{equation*}
i.e.,
\begin{equation*}
 _{x}D_{x_R}^{\alpha}u(x_i)=\delta'_{\alpha,_-x}{u_i^n} +\frac{1}{(\Delta x)^2}\epsilon(x_i)  +\mathcal{O}\left((\Delta
  x)^2\right).
\end{equation*}
 Next we calculate the error $\epsilon(x_i)$. Extending the
 definition of $u(x)$ from $[x_L,x_R]$ to $[x_L,x_R+\Delta x]$ with the extended part of $u(x)$ being the one given in Theorem A
 of Appendix, then
  \begin{equation*}
 \begin{split}
 \epsilon(x_i)= & \int_{x_{i-1}}^{x_R} \left(u(\xi)-S_{i-1}\left(\xi\right)\right)(\xi-x_{i-1})^{1-\alpha}d\xi \\
                &-2\int_{x_{i}}^{x_R} (u(\xi)-S_{i}(\xi))(\xi-x_{i})^{1-\alpha}d\xi
                +\int_{x_{i+1}}^{x_R} (u(\xi)-S_{i+1}(\xi))(\xi-x_{i+1})^{1-\alpha}d\xi   \\
              = & \int_{x_{i}}^{x_R+\Delta x} (u(\xi-\Delta x)-S_{i-1}(\xi-\Delta x))(\xi-x_{i})^{1-\alpha}d\xi \\
                &-2\int_{x_{i}}^{x_R} (u(\xi)-S_{i}(\xi))(\xi-x_{i})^{1-\alpha}d\xi\\
                &+\int_{x_{i}}^{x_R-\Delta x} (u(\xi+\Delta x)-S_{i+1}(\xi+\Delta x))(\xi-x_{i})^{1-\alpha}d\xi \\
              = &\sum_{k=i}^{N_x-1}\int_{x_{k}}^{x_{k+1}} (u(\xi-\Delta x)-S_{i-1}(\xi-\Delta x))(\xi-x_{i})^{1-\alpha}d\xi\\
                &-2\sum_{k=i}^{N_x-1}\int_{x_{k}}^{x_{k+1}} (u(\xi)-S_{i}(\xi))(\xi-x_{i})^{1-\alpha}d\xi\\
                &+\sum_{k=i}^{N_x-1}\int_{x_{k}}^{x_{k+1}} (u(\xi+\Delta x)-S_{i+1}(\xi+\Delta
                x))(\xi-x_{i})^{1-\alpha}d\xi+C(x_i),
 \end{split}
\end{equation*}
where
 \begin{equation*}
 \begin{split}
C(x_i)=\int_{x_R-\Delta x}^{x_R}\!\!\!\! (u(\xi)-S_{i-1}(\xi))(\xi-x_{i-1})^{1-\alpha}d\xi
      \!  -\int_{x_R}^{x_R+\Delta x}
       \!\!\!\!(u(\xi)-S_{i+1}(\xi))(\xi-x_{i+1})^{1-\alpha}d\xi,
\end{split}
\end{equation*}
and $C(x_i)=0$. Using Taylor series expansion and denoting
$\theta=(\xi-x_k)/\Delta x$, from Lemma 3.1 we get
 \begin{equation*}
 \begin{split}
\epsilon(x_i)=&\sum_{k=i}^{N_x-1}\int_{x_{k}}^{x_{k+1}} (\xi-x_{i})^{1-\alpha}\left[\frac{\theta(1-\theta)(\Delta x)^2}{2!}\left(u''(\xi-\Delta x)-2u''(\xi)
                                       + u''(\xi+\Delta x)\right)\right.\\
              &\left.+\frac{\theta(1-\theta)(1-2\theta)(\Delta x)^3}{3!}\left(u'''(\xi-\Delta x)-2u'''(\xi)+ u'''(\xi+\Delta x)\right ) +\mathcal{O}
                                     \left((\Delta x)^4 \right) \right]d\xi.\\
 \end{split}
\end{equation*}
Further using the first mean value theorem for integration and
Taylor series expansion, there exist ${\tilde\eta}{^k}$,
${\tilde\eta}{_1^k}$, ${\tilde\eta}{_2^k}$, ${\tilde\eta}{_3^k}$,
and ${\tilde\eta}{_4^k}$ such that
 \begin{equation*}
 \begin{split}
 &|\epsilon(x_i)|\\
=&\left|\sum_{k=i}^{N_x-1} \left\{ \left[\frac{\theta(1-\theta)(\Delta  x)^2}{2!}\big ( u''(\tilde{\eta}^k-\Delta x)-2u''(\tilde{\eta}^k) \right.\right.
                                     + u''(\tilde{\eta}^k+\Delta x)\big )+\frac{\theta(1-\theta)(1-2\theta)}{3!}\right. \\
 &\left.\times(\Delta x)^3\left(u'''(\tilde{\eta}^k-\Delta x)-2u'''(\tilde{\eta}^k)+ u'''(\tilde{\eta}^k+\Delta x)\right ) +\mathcal{O}
                                     \left.\left((\Delta x)^4 \right) \right]\left.\int_{x_{k}}^{x_{k+1}}\!\!(\xi-x_{i})^{1-\alpha}d\xi\right\}\right|\\
=&\left|\sum_{k=i}^{N_x-1} \left\{\left[\frac{\theta(1-\theta)(\Delta x)^2}{2!}\times\frac{(\Delta x)^2}{2!}
                     (u''''({\tilde\eta}{_1^k})\right.+u''''({\tilde\eta}{_2^k}))+\frac{\theta(1-\theta)(1-2\theta)(\Delta x)^4}{3!}\right.\right.\\
 &\left.\times(u''''({\tilde\eta}{_3^k})-u''''({\tilde\eta}{_4^k})) +\mathcal{O}
                                     \left.\left((\Delta x)^4 \right) \right]\left.\int_{x_{k}}^{x_{k+1}}(\xi-x_{i})^{1-\alpha}d\xi\right\}\right|\\
\leq  & |u''''({\tilde\eta}{_0})|\left[\frac{\theta(1-\theta)(\Delta x)^4}{2}\!+\frac{\theta(1-\theta)(1-2\theta)(\Delta x)^4}{3} \!+\mathcal{O}
                                     \left((\Delta x)^4 \right) \right]\sum_{k=i}^{N_x-1}\!\!\int_{x_{k}}^{x_{k+1}}\!\!\!(\xi-x_{i})^{1-\alpha}d\xi \\
=& |u''''({\tilde\eta}{_0})|\cdot\left[\frac{\theta(1-\theta)(\Delta x)^4}{2}+\frac{\theta(1-\theta)(1-2\theta)(\Delta x)^4}{3} +\mathcal{O}
                                     \left((\Delta x)^4 \right) \right]\frac{1}{2-\alpha}(x_R-x_i)^{2-\alpha},
 \end{split}
\end{equation*}
where $|u''''({\tilde\eta}{_0})|= \max \limits _{i\leq k \leq
N_x-1,\,1\leq j \leq 4} |u''''({\tilde\eta}{_j^k})|$.
Therefore $\epsilon(x_i)/(\Delta x)^2=\mathcal{O}\left((\Delta
x)^2\right)$, the desired result is proved.

\end{proof}

\noindent{\bf Remark 3.1.} From Lemma 3.3, it can be noted that $u\in C^4(\bar{\Omega})$ and  $_xD_{x_R}^{\alpha-2}u\in C^4(\bar{\Omega})$
 if and only if $u\in C^4(\bar{\Omega})$ and  $u^{(k)}(x_R)=0$, $k=0,1,2,3$. At the same time, if we do the zero extension to $u$ at $[x_R,x_R+\Delta x]$, then $u \in C^3[x_L,x_R+\Delta x]$ not belongs to $C^4[x_L,x_R+\Delta x]$.

Using similar idea we can prove

\noindent{\bf Lemma 3.5.}
  Let $u\in C^4(\bar{\Omega})$ and  $_{x_L}D_{x}^{\alpha-2}u\in C^4(\bar{\Omega})$,
  then  $\displaystyle_{x_L}D_{x}^{\alpha}u(x_i,t_n)= \delta'_{\alpha,_+x}{u_i^n}+\mathcal{O}\left((\Delta x)^2\right)$, where $\alpha \in (1,2)$.

\noindent{\bf Remark 3.2.} From Lemma 3.4, it can be noted that $u\in C^4(\bar{\Omega})$ and  $_{x_L}D_{x}^{\alpha-2}u\in C^4(\bar{\Omega})$
 if and only if $u\in C^4(\bar{\Omega})$ and  $u^{(k)}(x_L)=0$, $k=0,1,2,3$.

\noindent{\bf Example 1.} To numerically verify the truncation error given in Theorem 3.4 in a bounded domain, consider the function $u(x)=sin((1-x)^4)$, $x \in \Omega=(0,1)$; and it is easy to check that $u^{(4)}(x)|_{\partial \Omega} \neq 0$. Denoting $f(x)={_{x}D}_{x_R}^{\alpha}u(x)$, by the algorithm given in \cite{Deng:07}, we can numerically obtain the value of $f(x)$
at anywhere of the considered rectangle domain with any desired accuracy.

\begin{table}[h]\fontsize{9.5pt}{12pt}\selectfont
  \begin{center}
  \caption {The maximum errors (\ref{4.47}) and convergent orders for the scheme (\ref{2.18}).} \vspace{5pt}

    \begin{tabular*}{\linewidth}{@{\extracolsep{\fill}}*{7}{c}}                                    \hline  
         $\Delta x$ & $\alpha=1.1$ & Rate   & $\alpha=1.5$& Rate   & $\alpha=1.9$&   Rate      \\\hline
              ~~1/50&   1.0303e-002&        &  2.4122e-002&        &  3.8358e-002&      \\\hline 
             ~~1/100&   2.7832e-003&  1.8882&  6.4914e-003&  1.8937&  1.0322e-002&  1.8938    \\\hline 
             ~~1/200&   7.2217e-004&  1.9464&  1.6846e-003&  1.9461&  2.7038e-003&   1.9327    \\\hline 
             ~~1/400&   1.8379e-004&  1.9743&  4.2937e-004&  1.9721&  6.9889e-004&  1.9518    \\\hline 
    \end{tabular*}\label{tab:a}
  \end{center}
\end{table}

%

\noindent{\bf Example 2.} To numerically verify the truncation error given in Lemma 3.5 in a bounded domain, consider the function $u(x)=sin(x^4)$, $x \in \Omega=(0,1)$; and it is easy to check that $u^{(4)}(x)|_{\partial \Omega} \neq 0$. Denoting $F(x)={_{x_L}D}_{x}^{\alpha}u(x)$, again by the algorithm given in \cite{Deng:07}, we can numerically obtain the value of $F(x)$
at anywhere of the considered rectangle domain with any desired accuracy.

\begin{table}[h]\fontsize{9.5pt}{12pt}\selectfont
  \begin{center}
  \caption {The maximum errors (\ref{4.47}) and convergent orders for the scheme (\ref{2.19}). } \vspace{5pt}

    \begin{tabular*}{\linewidth}{@{\extracolsep{\fill}}*{7}{c}}                                    \hline  
         $\Delta x$ & $\alpha=1.1$ & Rate   & $\alpha=1.5$& Rate   & $\alpha=1.9$&   Rate      \\\hline
              ~~1/50&    1.0303e-002&        &  2.4122e-002&        &   3.8358e-002&             \\\hline 
             ~~1/100&    2.7832e-003&  1.8882&  6.4914e-003&  1.8937&  1.0322e-002&   1.8938    \\\hline 
             ~~1/200&   7.2217e-004&  1.9464&  1.6846e-003&  1.9461&   2.7038e-003&  1.9327    \\\hline 
             ~~1/400&    1.8379e-004&  1.9743&  4.2937e-004&  1.9721&  6.9889e-004&   1.9518    \\\hline 
    \end{tabular*}\label{tab:a}
  \end{center}
\end{table}

Table 1 and Table 2 numerically verify Theorem 3.4 and Lemma 3.5, respectively, and they show that the truncation errors are second order. In fact, watch carefully the numerical results, which also confirm Remark 2.1 in some sense.



\noindent{\bf Lemma 3.6.}
The coefficients $q_{i,k}^{\alpha}$ defined in (\ref{2.17}) satisfy
$\sum\limits_{k=i-1}^{N_x}q_{i,k}^{\alpha}< 0$,  $1 \leq i \leq N_x-1.$
\begin{proof}
Taking $u(x,t) \equiv 1$,  from (\ref{1.11}, \ref{2.155}, \ref{2.16}),
we have
$$\sum\limits_{k=i-1}^{N_x}q_{i,k}^{\alpha}=\sum\limits_{k=i-1}^{N_x}u(x_k,t)q_{i,k}^{\alpha}
=I_{\alpha}(x_{i-1})-2I_{\alpha}(x_i)+I_{\alpha}(x_{i+1}),$$
since $u(x,t) \equiv 1$ implies its linear interpolation function
$S_i(x)=u(x,t)\equiv 1$.

Then
\begin{equation*}
\begin{split}
       &I_{\alpha}(x_{i-1})-2I_{\alpha}(x_i)+I_{\alpha}(x_{i+1})\\
       &\quad=\frac{1}{\Gamma (2-\alpha)}\left[\int_{x_{i-1}}^{x_R}(\xi-{x_{i-1}})^{1-\alpha}d\xi
                -2\int_{x_i}^{x_R}(\xi-x_i)^{1-\alpha}d\xi+\int_{x_{i+1}}^{x_R}(\xi-{x_{i+1}})^{1-\alpha}d\xi \right]\\
       &\quad=\frac{1}{\Gamma (3-\alpha)}\left((x_R-x_{i-1})^{2-\alpha}
                -2(x_R-x_{i})^{2-\alpha}+(x_R-x_{i+1})^{2-\alpha} \right)  \\
       &\quad=\frac{1}{\Gamma (3-\alpha)}\left((a+\Delta x)^{2-\alpha}
                -2a^{2-\alpha}+(a-\Delta x)^{2-\alpha} \right)  \\
       &\quad=\frac{1}{\Gamma (3-\alpha)}\left[\left(\left(a+\Delta x\right)^{2-\alpha}-a^{2-\alpha}\right)
                -\left(a^{2-\alpha}-(a-\Delta x)^{2-\alpha} \right)\right]  \\
       &\quad=\frac{1}{(2-\alpha)\Gamma (3-\alpha)}\left(\int_a^{a+\Delta x} x^{1-\alpha}dx-\int_{a-\Delta x}^{a} x^{1-\alpha}dx\right)
\end{split}
\end{equation*}
where $a=x_R-x_i>0$; obviously, $x^{1-\alpha}$ is a decreasing
function, so
$I_{\alpha}(x_{i-1})-2I_{\alpha}(x_i)+I_{\alpha}(x_{i+1})<0.$
\end{proof}
Using similar method, we can prove

\noindent{\bf Lemma 3.7.}
The coefficients $p_{i,k}^{\alpha}$ defined in (\ref{2.20})  satisfy
$\sum\limits_{k=0}^{i+1}p_{i,k}^{\alpha}< 0$,   $1 \leq i \leq N_x-1.$

\noindent{\bf Theorem 3.8.}
The scheme (\ref{2.26}) of the fractional convection diffusion
equation (\ref{2.21}) with constant coefficients and $1<\alpha<2$ is unconditionally  von Neumann  stable.

\begin{proof}
Let $\widetilde{u_j^n}$ be the approximate solution of $u_j^n$,
which is the exact solution of the scheme (\ref{2.26}). Setting
$\varepsilon _j^n=\widetilde{u_j^n}- u_j^n$, $1 \leq j \leq N_x-1$,  then from (\ref{2.26})
we get the following perturbation equation
\begin{equation}\label{3.41}
 \left (1-\frac{\Delta t}{2} \delta_{\alpha,x}   \right )\varepsilon_j^{n+1}=
 \left (1+\frac{\Delta t}{2} \delta_{\alpha,x} \right )\varepsilon_j^{n},
\end{equation}
with zero boundary conditions, i.e.,
$\varepsilon_0^n=\varepsilon_{N_x}^n=0$. Then we can use the Von Neumann analysis
 or  Fourier
method \cite{Gustafsson:95} to do the stability analysis. Putting
$\xi_j=\frac{\Delta t}{2\Gamma\left(4-\alpha\right)\Delta
x^{\alpha}}d_{+,j}$, $\eta_j=\frac{\Delta t}{2\Gamma(4-\alpha)\Delta
x^{\alpha}}d_{-,j}$, $\gamma_j=\frac{\Delta t}{4\Delta x}g_j$,
 and
assuming
\begin{equation*}\label{3.42}
  \varepsilon_j^n=\frac{1}{\sqrt{2\pi}}\widehat{\varepsilon}^n(\omega)e^{i\omega
  x_j},
\end{equation*}
then (\ref{3.41}) leads to
\begin{equation*}
\begin{split}
&\widehat{\varepsilon}^{n+1}(\omega)
\left[ 1- \left(\xi_j\sum_{k=0}^{j+1}p_{j,k}^{\alpha}e^{i \omega (k-j) \Delta x}
+\eta_j\sum_{k=j-1}^{N_x}q_{j,k}^{\alpha}e^{i \omega (k-j)\Delta x}
+\gamma_j \left(e^{i \omega \Delta x}-e^{-i  \omega\Delta x} \right ) \right ) \right]\\
&\quad=\widehat{\varepsilon}^{n}(\omega)
\left[ 1+\left(\xi_j\sum_{k=0}^{j+1}p_{j,k}^{\alpha}e^{i \omega (k-j) \Delta x}
+\eta_j\sum_{k=j-1}^{N_x}q_{j,k}^{\alpha}e^{i \omega (k-j)\Delta x}
+\gamma_j \left(e^{i \omega \Delta x}-e^{-i  \omega\Delta x} \right ) \right ) \right].\\
\end{split}
\end{equation*}
The amplification factor is
\begin{equation*}
\begin{split}
&\widehat{Q}(\omega)=\left[ 1- \left(\xi_j\sum_{k=0}^{j+1}p_{j,k}^{\alpha}e^{i \omega (k-j) \Delta x}
+\eta_j\sum_{k=j-1}^{N_x}q_{j,k}^{\alpha}e^{i \omega (k-j)\Delta x}
+\gamma_j \left(e^{i \omega \Delta x}-e^{-i  \omega\Delta x} \right ) \right ) \right]^{-1}\\
&\qquad \cdot \left[ 1+\left(\xi_j\sum_{k=0}^{j+1}p_{j,k}^{\alpha}e^{i \omega (k-j) \Delta x}
+\eta_j\sum_{k=j-1}^{N_x}q_{j,k}^{\alpha}e^{i \omega (k-j)\Delta x}
+\gamma_j \left(e^{i \omega \Delta x}-e^{-i  \omega\Delta x} \right ) \right ) \right].\\
\end{split}
\end{equation*}
Next we prove that  $|\widehat{Q}(\omega)|< 1$,  it means that  the real part of the element in the big
parentheses of above equation is negative, that is
\begin{equation*}
\sum_{k=0}^{j+1}p_{j,k}^{\alpha}cos( \omega (k-j) \Delta x) <0 ~\mbox {and} ~
\sum_{k=j-1}^{N_x}q_{j,k}^{\alpha}cos( \omega (k-j) \Delta x) < 0.
 \end{equation*}
We can write
\begin{equation*}
\begin{split}
\sum_{k=0}^{j+1}p_{j,k}^{\alpha}cos( \omega (k-j) \Delta x)
&=(p_{j,j-1}^{\alpha}+p_{j,j+1}^{\alpha})cos( \omega \Delta x)+p_{j,j}^{\alpha} +\sum_{k=0}^{j-2}p_{j,k}^{\alpha}cos( \omega (k-j) \Delta x),\\
\sum_{k=j-1}^{N_x}\!\!q_{j,k}^{\alpha}cos( \omega (k-j) \Delta x)&=(q_{j,j-1}^{\alpha}+q_{j,j+1}^{\alpha})cos( \omega \Delta
x)+q_{j,j}^{\alpha}+\sum_{k=j+2}^{N_x}q_{j,k}^{\alpha}cos( \omega
(k-j) \Delta x).
\end{split}
 \end{equation*}
From Theorem B of Appendix, we know
$p_{j,j-1}^{\alpha}+p_{j,j+1}^{\alpha}\geq 0$, $p_{j,k}^{\alpha}
\geq 0$ for $k \leq j-2$, and
$q_{j,j-1}^{\alpha}+q_{j,j+1}^{\alpha}\geq 0$,  $q_{j,k}^{\alpha}
\geq 0$ for $k \geq j+2$, then
\begin{equation*}
\begin{split}
\sum_{k=0}^{j+1}  p_{j,k}^{\alpha}cos( \omega (k-j) \Delta x)  &\leq (p_{j,j-1}^{\alpha}+p_{j,j+1}^{\alpha})+p_{j,j}^{\alpha}
+\sum_{k=0}^{j-2}p_{j,k}^{\alpha} \leq \sum_{k=0}^{j+1}p_{j,k}^{\alpha} <0,\\
\sum_{k=j-1}^{N_x}q_{j,k}^{\alpha}cos( \omega (k-j) \Delta x)  & \leq (q_{j,j-1}^{\alpha}+q_{j,j+1}^{\alpha}) +q_{j,j}^{\alpha}+\!\!\sum_{k=j+2}^{N_x}q_{j,k}^{\alpha}
\leq \!\! \sum_{k=j-1}^{N_x}q_{j,k}^{\alpha} < 0,
\end{split}
 \end{equation*}
where Lemma 3.6 and 3.7 are used.
\end{proof}

\subsection{Convergence and stability for two-dimensional two-sided fractional convection diffusion equation}
In this subsection, we prove the consistency and stability of the
scheme (\ref{2.33})-(\ref{2.37}) which combine the alternating
directions implicit scheme with Crank-Nicolson scheme together
(ADI-CN). According to the Lax's equivalence theorem, the
convergence of the ADI-CN scheme is naturally obtained.

\noindent{\bf Theorem 3.9.}
The truncation error of the scheme (\ref{2.33}) is
$\mathcal{O}((\Delta x)^2)+\mathcal{O}((\Delta
y)^2)+\mathcal{O}((\Delta t)^2)$.

\begin{proof}
  Let $u(x,y,t)$ be the exact solution of the  two-dimensional two-sided fractional convection diffusion equation (\ref{1.1}), then for the
 scheme (\ref{2.32}) we have
  \begin{equation*}
  \begin{split}
  &\frac{u_{i,j}^{n+1}-u_{i,j}^{n} }{\Delta t} -\big(d_{+,i,j}\delta'_{\alpha,_+x} + d_{-,i,j} \delta'_{\alpha,_-x}+e_{+,i,j}\delta'_{\beta,_+y}
                      + e_{-,i,j} \delta'_{\beta,_-y} \\
                      & + g_{i,j}D'_{\alpha,x}+ h_{i,j}D'_{\beta,y} \big )\frac{{u_{i,j}^{n+1}+u_{i,j}^{n}}}{2}- s_{i,j}^{n+1/2}\\
= &\frac{\partial u(x_i,y_j,t_{n+1/2})}{\partial t}
                      -d_{+,i,j}~ {_{x_L}D_{x}^{\alpha}u(x_i,y_j,t_{n+1/2})}
                      -d_{-,i,j}~ {_{x}D_{x_R}^{\alpha}u(x_i,y_j,t_{n+1/2})}\\
  & -e_{+,i,j}~ {_{y_L}D_{y}^{\beta}u(x_i,y_j,t_{n+1/2})}-e_{-,i,j}~ {_{y}D_{y_R}^{\beta}u(x_i,y_j,t_{n+1/2})}
                      -g_{i,j}u_x(x_i,y_j,t_{n+1/2})\\
     &     -h_{i,j}u_y(x_i,y_j,t_{n+1/2})-s(x_i,y_j,t_{n+1/2})           +\mathcal{O}(\Delta x)^2+\mathcal{O}(\Delta y)^2+\mathcal{O}(\Delta t)^2.
\end{split}
\end{equation*}
The scheme (\ref{2.33}) differs from (\ref{2.32}) by a perturbation
equals to \cite{Charles:07}
\begin{equation*}
  \frac{(\Delta
  t)^2}{4}\delta_{\alpha,x}\delta_{\beta,y}(u_{i,j}^{n+1}-u_{i,j}^n),
\end{equation*}
which may be deduced by distributing the operator products in
(\ref{2.33}). Since $(u_{i,j}^{n+1}-u_{i,j}^n)$ is an
$\mathcal{O}(\Delta t)$ term, it follows that the perturbation
contributes an $\mathcal{O}((\Delta t)^2)$ error component to the
truncation error of (\ref{2.32}). Thus, the scheme (\ref{2.33}) has
a truncation error also $\mathcal{O}((\Delta
x)^2)+\mathcal{O}((\Delta y)^2)+\mathcal{O}((\Delta t)^2)$.
\end{proof}

\noindent{\bf Theorem 3.10.}
  The ADI-CN scheme, defined by (\ref{2.32}) with constant coefficients, is unconditionally von Neumann  stable for   $1<\alpha,\beta <2$.

\begin{proof}
Let $\widetilde{u_{j,m}^n}$ be the approximate solution of $u_{j,m}^n$,
which is the exact solution of the scheme (\ref{2.32}). Setting
$\varepsilon _{j,m}^n=\widetilde{u_{j,m}^n}- u_{j,m}^n$, $1 \leq j \leq N_x-1$, and $1 \leq m \leq N_y-1$,
 then from (\ref{2.32})
we get the following perturbation equation
\begin{equation}\label{0036}
 \left (1-\frac{\Delta t}{2}\delta_{\alpha,x}-\frac{\Delta t}{2}\delta_{\beta,y}\right )\varepsilon_{j,m}^{n+1}
 =\left (1+\frac{\Delta t}{2}\delta_{\alpha,x}+\frac{\Delta t}{2}\delta_{\beta,y}\right ) \varepsilon_{j,m}^{n},
\end{equation}
Similarily,  putting
$\xi_j=\frac{\Delta t}{2\Gamma\left(4-\alpha\right)\Delta
x^{\alpha}}d_{+,j}$, $\eta_j=\frac{\Delta t}{2\Gamma(4-\alpha)\Delta
x^{\alpha}}d_{-,j}$, $\gamma_j=\frac{\Delta t}{4\Delta x}g_j$  and
$\tilde{\xi}_m=\frac{\Delta t}{2\Gamma\left(4-\beta\right)\Delta y^{\beta}}e_{+,m}$,
$\tilde{\eta_m}=\frac{\Delta t}{2\Gamma(4-\beta)\Delta
y^{\beta}}e_{-,m}$, $\tilde{\gamma_m}=\frac{\Delta t}{4\Delta y}h_m$,
assuming
\begin{equation*}
  \varepsilon_{j,m}^n=\frac{1}{\sqrt{2\pi}}\tilde{\varepsilon}^n(\omega)e^{i\omega (x_j+y_m)},
\end{equation*}
then (\ref{0036}) leads to
\begin{equation*}
\begin{split}
&\tilde{\varepsilon}^{n+1}(\omega)
\bigg[ 1- \left(\xi_j\sum_{k=0}^{j+1}p_{j,k}^{\alpha}e^{i \omega (k-j) \Delta x}
+\eta_j\sum_{k=j-1}^{N_x}q_{j,k}^{\alpha}e^{i \omega (k-j)\Delta x}
+\gamma_j \left(e^{i \omega \Delta x}-e^{-i  \omega\Delta x} \right ) \right )\\
&\qquad\qquad~\, - \left(\tilde{\xi}_m\sum_{k=0}^{m+1}p_{m,k}^{\beta}e^{i \omega (k-m) \Delta y}
+\tilde{\eta_m}\sum_{k=m-1}^{N_y}q_{m,k}^{\beta}e^{i \omega (k-m)\Delta y}
+\tilde{\gamma_m} \left(e^{i \omega \Delta y}-e^{-i  \omega\Delta y} \right ) \right ) \bigg]\\
&\quad=\tilde{\varepsilon}^{n}(\omega)
\bigg[ 1+ \left(\xi_j\sum_{k=0}^{j+1}p_{j,k}^{\alpha}e^{i \omega (k-j) \Delta x}
+\eta_j\sum_{k=j-1}^{N_x}q_{j,k}^{\alpha}e^{i \omega (k-j)\Delta x}
+\gamma_j \left(e^{i \omega \Delta x}-e^{-i  \omega\Delta x} \right ) \right )\\
&\qquad\qquad~\, + \left(\tilde{\xi}_m\sum_{k=0}^{m+1}p_{m,k}^{\beta}e^{i \omega (k-m) \Delta y}
+\tilde{\eta_m}\sum_{k=m-1}^{N_y}q_{m,k}^{\beta}e^{i \omega (k-m)\Delta y}
+\tilde{\gamma_m} \left(e^{i \omega \Delta y}-e^{-i  \omega\Delta y} \right ) \right ) \bigg].\\
\end{split}
\end{equation*}
The amplification factor is
\begin{equation*}
\begin{split}
&\tilde{Q}(\omega)=\bigg[ 1- \left(\xi_j\sum_{k=0}^{j+1}p_{j,k}^{\alpha}e^{i \omega (k-j) \Delta x}
+\eta_j\sum_{k=j-1}^{N_x}q_{j,k}^{\alpha}e^{i \omega (k-j)\Delta x}
+\gamma_j \left(e^{i \omega \Delta x}-e^{-i  \omega\Delta x} \right ) \right )\\
&\qquad - \left(\tilde{\xi}_m\sum_{k=0}^{m+1}p_{m,k}^{\beta}e^{i \omega (k-m) \Delta y}
+\tilde{\eta_m}\sum_{k=m-1}^{N_y}q_{m,k}^{\beta}e^{i \omega (k-m)\Delta y}
+\tilde{\gamma_m} \left(e^{i \omega \Delta y}-e^{-i  \omega\Delta y} \right ) \right ) \bigg]^{-1}\\
&\qquad \cdot \bigg[ 1+ \left(\xi_j\sum_{k=0}^{j+1}p_{j,k}^{\alpha}e^{i \omega (k-j) \Delta x}
+\eta_j\sum_{k=j-1}^{N_x}q_{j,k}^{\alpha}e^{i \omega (k-j)\Delta x}
+\gamma_j \left(e^{i \omega \Delta x}-e^{-i  \omega\Delta x} \right ) \right )\\
&\qquad + \left(\tilde{\xi}_m\sum_{k=0}^{m+1}p_{m,k}^{\beta}e^{i \omega (k-m) \Delta y}
+\tilde{\eta_m}\sum_{k=m-1}^{N_y}q_{m,k}^{\beta}e^{i \omega (k-m)\Delta y}
+\tilde{\gamma_m} \left(e^{i \omega \Delta y}-e^{-i  \omega\Delta y} \right ) \right ) \bigg].\\
\end{split}
\end{equation*}
Next we prove that $|\widehat{Q}(\omega)|< 1$,  it means that  the real part of the element in the big
parentheses of above equation is negative, that is
\begin{equation*}
\begin{split}
\sum_{k=0}^{j+1}p_{j,k}^{\alpha}cos( \omega (k-j) \Delta x) <0 ~ & \mbox {and} ~
\sum_{k=j-1}^{N_x}q_{j,k}^{\alpha}cos( \omega (k-j) \Delta x) < 0,\\
\sum_{k=0}^{m+1}p_{m,k}^{\beta}cos( \omega (k-m) \Delta y) <0 ~  &\mbox {and} ~
\sum_{k=m-1}^{N_y}q_{m,k}^{\beta}cos( \omega (k-m) \Delta y) < 0.
\end{split}
 \end{equation*}
Similar to Theorem  3.8, we get that the scheme
(\ref{2.32}) is unconditionally von Neumann stable.
\end{proof}

 \section{Numerical Results}
Here we verify the above theoretical results including convergent
order and stability.  Introducing
the vectors $U=[u_h(x_0,t),\ldots,u_h(x_N,t)]^{\rm T}$, where $U$ is
the approximated value,  and $u=[u(x_0,t),\ldots,u(x_N,t)]^{\rm T}$,
where $u$ is the exact value, in the following numerical examples
the errors are measured by
\begin{equation}\label{4.47}
  ||U(\Delta x)-u(\Delta x)||_\infty,
\end{equation}
where $||\cdot||_\infty$ is the $ l_\infty$ norm.

\subsection{Numerical results for one-dimensional two-sided fractional convection diffusion equation}
Let us consider  the one-dimensional fractional convection diffusion
equation (\ref{2.21}),
where $0< x < 1 $ and $0 < t \leq1$, with the coefficient functions
\begin{equation*}
\begin{split}
  &d_+(x)=\Gamma(3-\alpha)x^{\alpha},\quad d_-(x) =\Gamma(3-\alpha)(2-x)^{\alpha},\quad\mbox{and} \quad g(x)=\frac{1}{4}x.
\end{split}
\end{equation*}
Take the exact solution of the equation as $u(x,t)=e^{-t}sin((2x)^4)sin((2-2x)^4)$, then the corresponding initial and boundary conditions are, respectively, $u(x,0)=sin((2x)^4)sin((2-2x)^4)$ and $u(0,t)=u(1,t)=0$; and the forcing function
\begin{equation*}
\begin{split}
 s(x,t)=&-e^{-t}sin((2x)^4)sin((2-2x)^4)-e^{-t}d_{+}(x)\,{_{x_L}D_x^{\alpha}sin((2x)^4)sin((2-2x)^4)}\\
 &-d_{-}(x)\,{_{x}D_{x_R}^{\alpha}sin((2x)^4)sin((2-2x)^4)}\\
 &-64e^{-t}g(x)(x^3cos((2x)^4)sin( (2-2x)^4)  -(1-x)^3sin((2x)^4)cos( (2-2x)^4)  ),
  \end{split}
\end{equation*}
by the algorithm given in \cite{Deng:07}, we can numerically obtain the value of $s(x,t)$
at anywhere of the considered rectangle domain with any desired accuracy.

\begin{table}[h]\fontsize{9.5pt}{12pt}\selectfont
  \begin{center}
  \caption {The maximum errors (\ref{4.47}) and convergent orders for the scheme (\ref{2.22}) of the one-dimensional two-side fractional convection diffusion equation (\ref{2.21}) at t=1 and $\Delta t=\Delta x$.} \vspace{5pt}

    \begin{tabular*}{\linewidth}{@{\extracolsep{\fill}}*{7}{c}}                                    \hline  
         $\Delta t,\,\Delta x$ & $\alpha=1.1$ & Rate   & $\alpha=1.5$ & Rate   & $\alpha=1.9$ &   Rate      \\\hline
              ~~1/50&   2.1180e-003&        &  1.9815e-003&        &  1.3809e-003&             \\\hline 
             ~~1/100&   5.2688e-004&  2.0072&  5.0092e-004&  1.9839&  3.5593e-004&   1.9559    \\\hline 
             ~~1/200&   1.3174e-004&  1.9997&  1.2649e-004&  1.9856&  9.1681e-005&   1.9569    \\\hline 
             ~~1/400&   3.2913e-005&  2.0010&  3.1851e-005&  1.9896&  2.3523e-005&   1.9625    \\\hline 
    \end{tabular*}\label{tab:a}
  \end{center}
\end{table}

In Table 3, we show the scheme (\ref{2.22}) is second order
convergent in both space and time directions.

\subsection{Numerical results for two-dimensional two-sided fractional convection diffusion equation}
 Consider the two-dimensional two-sided fractional convection diffusion equation (\ref{1.1}),
on a finite domain $ 0<  x< 1,\,0<  y< 1$, $0< t \leq 1$, and with the
coefficients
 \begin{equation*}
\begin{split}
&d_+(x,y)=\Gamma(3-\alpha)x^{\alpha},\quad d_-(x,y)=\Gamma(3-\alpha)(2-x)^{\alpha},\quad g(x,y)=\frac{1}{4}x,  \\
&e_+(x,y)=\Gamma(3-\beta)y^{\beta},\quad e_-(x,y)=\Gamma(3-\beta)(2-y)^{\beta}, \,\quad h(x,y)=\frac{1}{4}y,   \\
\end{split}
\end{equation*}
and the initial condition $u(x,y,0)=sin((2x)^4)sin((2-2x)^4)sin((2y)^2)sin((2-2y)^2)$ and the
Dirichlet boundary conditions on the rectangle in the form $
u(0,y,t)=u(x,0,t)=0$ and $ u(1,y,t)=u(x,1,t)=0$
 for all $ t\geq 0$. The exact solution to this two-dimensional two-sided fractional convection diffusion equation is
 \begin{equation*}
\begin{split}
u(x,y,t)=e^{-t}sin((2x)^4)sin((2-2x)^4)sin((2y)^4)sin((2-2y)^4).
  \end{split}
  \end{equation*}
By the algorithm given in \cite{Deng:07} and above conditions, it is easy to obtain the forcing function $s(x,y,t)$
at anywhere of the considered rectangle domain with any desired accuracy.

\begin{table}[h]\fontsize{9.5pt}{12pt}\selectfont
  \begin{center}
  \caption{The maximum errors (\ref{4.47}) and convergent orders for the scheme (\ref{2.33})-(\ref{2.37}) of
  the two-dimensional two-sided fractional convection diffusion equation (\ref{1.1}) at $t=1$ and $\Delta t
=\Delta x =\Delta y$.}\vspace{5pt}

    \begin{tabular*}{\linewidth}{@{\extracolsep{\fill}}*{7}{c}}                                    \hline  
         $\Delta t,\,\Delta x,\,\Delta y$ & $\alpha=1.1,\beta=1.1$ & Rate   & $\alpha=1.6,\beta=1.4$ & Rate   & $\alpha=1.9,\beta=1.9$ &   Rate      \\\hline
             ~~1/25&   9.5946e-003&        &  8.5313e-003&        &  1.0232e-002&             \\\hline 
             ~~1/50&   2.3956e-003&  2.0018&  2.1729e-003&  1.9731&  2.6207e-003&   1.9650    \\\hline 
             ~~1/100&  5.9582e-004&  2.0075&  5.5244e-004&  1.9757&  6.6155e-004&   1.9860    \\\hline 
             ~~1/200&  1.4915e-004&  1.9981&  1.3959e-004&  1.9847&  1.6796e-004&   1.9778    \\\hline 
    \end{tabular*}\label{tab:a}
  \end{center}
\end{table}

Table 4 also shows the  maximum error, at  time $t=1$ and $\Delta t
=\Delta x =\Delta y$, between the exact analytical value and the
numerical value obtained by applying the ADI-CN scheme
(\ref{2.33})-(\ref{2.37}), and the scheme is second order convergent
and this is in agreement with the order of the truncation error.

\subsection{Numerical results for two-dimensional one-sided fractional convection diffusion equation}
 Considering the two-dimensional two-sided fractional convection diffusion equation (\ref{1.1}), and taking the coefficients functions
 as
  \begin{equation*}
\begin{split}
&d_+(x,y)=1,\quad d_-(x,y)=0, \quad g(x,y)=1,\\
&e_+(x,y)=1, ~\quad
e_-(x,y)=0,\quad h(x,y)=1,
  \end{split}
\end{equation*}
then it becomes the two-dimensional one-sided fractional convection
diffusion equation,

\begin{equation}\label{4.49}
\begin{split}
\frac{\partial u(x,y,t) }{\partial t}=& \frac{\partial^{\alpha}u(x,y,t) }{\partial {x^{\alpha}}}+ \frac{\partial^{\beta}u(x,y,t) }{\partial {y^{\beta}}}
                                      + \frac{\partial u(x,y,t) }{\partial x}+\frac{\partial u(x,y,t) }{\partial
                                      y}+s(x,y,t),
\end{split}
\end{equation}
where $ 0<  x < 1$,  $0<  y < 1$, $0< t \leq 1$,
and the initial condition $u(x,y,0)=sin(x^4)sin(y^4)$ and the Dirichlet
boundary conditions on the rectangle in the simple  form $
u(0,y,t)=u(x,0,t)=0,u(1,y,t)=e^{-t}sin(1)sin(y^4)$ and $u(x,1,t)=e^{-t}sin(1)sin(x^4)$ for
all $ t> 0$. The exact value to this two-dimensional one-sided
fractional convection diffusion equation is
 \begin{equation*}
\begin{split}
u(x,y,t)=e^{-t}sin(x^4)sin(y^4).
  \end{split}
  \end{equation*}
By the algorithm given in \cite{Deng:07} and above conditions, it is easy to obtain the forcing function $s(x,y,t).$
at anywhere of the considered rectangle domain with any desired accuracy.

\begin{table}[h]\fontsize{9.5pt}{12pt}\selectfont
  \begin{center}
  \caption{The maximum errors (\ref{4.47}) and convergent orders for the scheme (\ref{2.33})-(\ref{2.37}) of the two-dimensional one-sided fractional convection diffusion equation (\ref{4.49}) at $t=1$ and $\Delta t
=\Delta x =\Delta y$.}\vspace{5pt}

    \begin{tabular*}{\linewidth}{@{\extracolsep{\fill}}*{7}{c}}                                    \hline  
         $\Delta t,\,\Delta x,\,\Delta y$ & $\alpha=1.1,\beta=1.1$ & Rate   & $\alpha=1.6,\beta=1.4$ & Rate   & $\alpha=1.9,\beta=1.9$ &   Rate      \\\hline
             ~~1/25&   1.1435e-003&        &  5.0896e-004&        &  2.6381e-004&             \\\hline 
             ~~1/50&   2.8953e-004&  1.9817&  1.3592e-004&  1.9048&  6.9390e-005&   1.9267    \\\hline 
             ~~1/100&  6.8091e-005&  2.0882&  3.4877e-005&  1.9624&  1.8064e-005&   1.9416    \\\hline 
             ~~1/200&  1.5950e-005&  2.0939&  8.8502e-006&  1.9785&  4.6728e-006&   1.9507    \\\hline 
    \end{tabular*}\label{tab:a} \vspace{20pt}
  \end{center}
\end{table}
Table 5 shows the maximum error, at  time $t=1$ and $\Delta t
=\Delta x =\Delta y$, between the exact analytical value and the
numerical value obtained by applying the ADI-CN scheme
(\ref{2.33})-(\ref{2.37}), and the scheme is second order convergent
and  it  corresponds to the order of the truncation error.

\section{Conclusions}
This work provides the second-order efficient numerical scheme for
the two-dimensional two-sided fractional advection diffusion
equation on a finite domain. Both the convergent order and numerical
stability of the scheme are theoretically proved and numerically
verified.
This paper can be considered as the sequel of the works
\cite{Meerschaert:04,Meerschaert:06,Mark:06,Sousa:Pro,Sousa:11,Sousa:12,Charles:07,Charles:06}.

\section*{Acknowledgments} This work was supported by the Program for
New Century Excellent Talents in University under Grant No.
NCET-09-0438, the National Natural Science Foundation of China under
Grant No. 11271173, and the Fundamental Research Funds for the
Central Universities under Grant No. lzujbky-2010-63 and lzujbky-2012-k26. M.H. Chen
thanks Meerschaert for providing the materials and discussions.

\section*{Appendix}
\newtheorem*{Append1}{Theorem A} \label{Theorem 1}

\begin{Append1}
{\em For any given $u(x) \in C^4[x_L,x_R]$, there exists an
extension of $u(x)$ defined on $[x_R,x_R+\Delta x]$, such that the
extended $u(x) \in C^4[x_L,x_R+\Delta x]$ and
\begin{equation} \label{appendix 1}
\int_{x_R-\Delta x}^{x_R}
(u(\xi)-S_{i-1}(\xi))(\xi-x_{i-1})^{1-\alpha}d\xi
       -\int_{x_R}^{x_R+\Delta x}
       (u(\xi)-S_{i+1}(\xi))(\xi-x_{i+1})^{1-\alpha}d\xi=0,
\end{equation}
where $S_{i-1}(x)$ is the linear interpolation function of $u(x)$ on
the interval $[x_R-\Delta x,x_R]$, and $S_{i+1}(x)$ the linear
interpolation function of the extended $u(x)$ on the interval $[x_R,x_R+\Delta
x]$.}
\end{Append1}
\begin{proof}
Denoting $a=x_R$, $b=x_R+\Delta x$, and $d_1 =\int_{x_R-\Delta
x}^{x_R} (u(\xi)-S_{i-1}(\xi))(\xi-x_{i-1})^{1-\alpha}d\xi$, then
the equality (\ref{appendix 1}) means to find the extension of
$u(x)$ on the interval $[a,b]$ such that $\int_{a}^{b}
(u(\xi)-S_{i+1}(\xi))(\xi-x_{i+1})^{1-\alpha}d\xi= d_1$. Taking
$u(b)=u(a)$, then $S_{i+1}(x) \equiv u(a)$ on the interval $[a,b]$
since $S_{i+1}(x)$ is the linear interpolation of $ u(x)$. Now we
need to prove that there exists $V(x)=u(x)-S_{i+1}(x)$  on the
interval $[a,b]$ such that $V(a)=0,$ $V(b)=0$,
$V'(a^+)=u'(a^-)=d_2$, $V''(a^+)=u''(a^-)=d_3$,
$V'''(a^+)=u'''(a^-)=d_4$, and  $V''''(a^+)=u''''(a^-)=d_5$.

Suppose that $V(x)$ is composed by
$\{1,x,x^2,x^2,x^3,x^4,x^5,x^6\}$, i.e.,
 \begin{equation*}
   V(x)=c_0+c_1x+c_2x^2+c_3x^3+c_4x^4+c_5x^5+c_6x^6.
   \end{equation*}
 Using the first mean value theorem for integration, there exists $\eta \in(a,b)$, such that
 \begin{equation*}
 \begin{split}
 \int_{a}^{b} V(\xi)(\xi-x_{i+1})^{1-\alpha}d\xi= V(\eta)\int_{a}^{b}
 (\xi-x_{i+1})^{1-\alpha}d\xi=d_1,
\end{split}
\end{equation*}
  and it can be rewritten as
\begin{equation*}
c_0+c_1\eta+c_2\eta^2+c_3\eta^3+c_4\eta^4+c_5\eta^5+c_6\eta^6=\left(\int_{a}^{b}
(\xi-x_{i+1})^{1-\alpha}d\xi \right)^{-1}d_1 =:d_1
\end{equation*}
 (again denoted
by $d_1$, and it is easy to check that $\int_{a}^{b}
(\xi-x_{i+1})^{1-\alpha}d\xi \neq 0$). Combining all the above
requirements, we obtain the linear algebraic equations $Ac=d$, where
$c=[c_0,c_1,c_2,c_3,c_4,c_5,c_6]^{\rm T}$,
$d=[0,0,d_1,d_2,d_3,d_4,d_5]^{\rm T}$, and the determinant of the
coefficient matrix $A$ is

\begin{equation*}
det(A)=\left| \begin{array}
 {l@{\qquad} c@{\qquad} c@{\qquad} c@{\qquad} c@{\qquad} c@{\qquad} c}
 1  & a  & a^2  & a^3    & a^4   &a^5    & a^6  \\
 1  & b  & b^2  & b^3    & b^4   &b^5    & b^6  \\
 1  & \eta  & \eta^2  & \eta^3    & \eta^4   &\eta^5    & \eta^6  \\
 0  & 1  & 2a   & 3a^2   & 4a^3  &5a^4   &6a^5  \\
 0  & 0  & 2    & 6a     &12a^2  &20a^3  &30a^4  \\
 0  & 0  & 0    & 6      &24a    &60a^2  &120a^3  \\
 0  & 0  & 0    & 0      &24     &120a   &360a^2  \\
\end{array}
 \right|
\end{equation*}
\begin{equation*}
~~~~~~~~=\left| \begin{array}
 {l@{\qquad} c@{\qquad} c@{\qquad} c@{\qquad} c@{\qquad} c@{\qquad} c}

 1  & a    & a^2      & a^3        & a^4       &a^5        & a^6  \\
 0  & b-a  & b^2-a^2  & b^3-a^3    & b^4-a^4   &b^5-a^5    & b^6-a^6  \\
 0  & \eta-a  & \eta^2-a^2  & \eta^3-a^3    & \eta^4-a^4   &\eta^5-a^5    & \eta^6-a^6  \\
 0  & 1    & 2a       & 3a^2       & 4a^3      &5a^4   &6a^5  \\
 0  & 0    & 2        & 6a         &12a^2      &20a^3  &30a^4  \\
 0  & 0    & 0        & 6          &24a        &60a^2  &120a^3  \\
 0  & 0    & 0        & 0          &24         &120a   &360a^2
\end{array}
 \right|.
\end{equation*}
Using the formula
$a^n-b^n=(a-b)(a^{n-1}+ba^{n-2}+b^2a^{n-3}+\ldots+b^{n-2}a+b^{n-1})$,
taking out the common factors $\eta-a$ and $b-a$, and repeating the
steps, finally we get
  \begin{equation*}
     \begin{split}
det(A)=&288(\eta-a)^5(b-a)^5\left| \begin{array} {l@{\qquad}
c@{\qquad} c@{\qquad} c@{\qquad} c@{\qquad} c@{\qquad} c}
 1  & a  & a^2  & a^3     & a^4    &a^5      & a^6  \\
 0  & 1  & 2a   & 3a^2    & 4a^3   &5a^4     & 6a^5  \\
 0  & 0  & 1    & 3a      & 6a^2   &10a^3    & 15a^4  \\
 0  & 0  & 0    & 1       & 4a     &10a^2    &20a^3  \\
 0  & 0  & 0    & 0       &1       &5a       &15a^2  \\
 0  & 0  & 0    & 0       &0       &1        &b+5a  \\
 0  & 0  & 0    & 0       &0       &0        &\eta-b  \\
 \end{array}
 \right|\\
 =&288(\eta-b)(\eta-a)^5(b-a)^5.
 \end{split}
\end{equation*}
Since $\eta\neq a$, $\eta\neq b$, and $a \neq b$, we get
$det(A)\neq0$. Then with the given basis functions $\{1$, $x$,
$x^2$, $x^2$, $x^3$, $x^4$, $x^5$, $x^6\}$, there exists a unique
extended part of $u(x)$, which is $V(x)+u(a)$ with $x\in[x_R,x_R+
\Delta x]$. The extended $u(x)\in C^4[x_L,x_R+ \Delta x]$ and
(\ref{appendix 1}) holds.
\end{proof}

\newtheorem*{Append2}{Theorem B}
\begin{Append2}
  {\em For the coefficients $q_{i,k}^{\alpha}$ defined in (\ref{2.17}) and $p_{i,k}^{\alpha}$ defined by
  (\ref{2.20}), the following hold:
   $q_{i,k}^{\alpha}> 0$ for $k \geq i+2$, $p_{i,k}^{\alpha}> 0$ for $k \leq i-2$, and $p_{i,i+1}^{\alpha} =q_{i,i-1}^{\alpha} =1$,
   $p_{i,i}^{\alpha}= q_{i,i}^{\alpha}=-4+2^{3-\alpha}$, $p_{i,i-1}^{\alpha} =q_{i,i+1}^{\alpha} =6-2^{5-\alpha}+3^{3-\alpha} $.}
\end{Append2}
\begin{proof}
Similar to the proof given in \cite{Sousa:12} for the left
Riemann-Liouville fractional derivative, we have
   $q_{i,k}^{\alpha}=0,\,\,k<i-1,\,\,
  q_{i,i-1}^{\alpha} =1,\,\,
  q_{i,i}^{\alpha}=-4+2^{3-\alpha},\,\,
  q_{i,i+1}^{\alpha}\ =6-2^{5-\alpha}+3^{3-\alpha}, \,\, \mbox{and}\,\,
  q_{i,k}^{\alpha}=(k-i+2)^{3-\alpha}-4(k-i+1)^{3-\alpha}+6(k-i)^{3-\alpha}-4(k-i-1)^{3-\alpha}+(k-i-2)^{3-\alpha},\,\,k\geq i+2$.
  For  $ q_{i,k}^{\alpha},k\geq i+2$, denote $ m=k-i\geq 2$, then
    \begin{equation*}
  \begin{split}
   q_{i,i+m}^{\alpha}&=(m+2)^{3-\alpha}-4(m+1)^{3-\alpha}+6m^{3-\alpha}-4(m-1)^{3-\alpha}+(m-2)^{3-\alpha}\\
                     &=m^{3-\alpha}\left[6+ \left(1+\frac{2}{m}\right)^{3-\alpha}  -4\left(1+\frac{1}{m}\right)^{3-\alpha}
                     -4\left(1-\frac{1}{m}\right)^{3-\alpha} +\left(1-\frac{2}{m}\right)^{3-\alpha}    \right]\\
                     &=m^{3-\alpha}\left\{6+ \sum_{k=0}^{\infty}{ 3-\alpha \choose k } \left[\left(\frac{2}{m}\right)^k  -4\left(\frac{1}{m}\right)^k -4\left(\frac{-1}{m}\right)^k +\left(\frac{-2}{m}\right)^k   \right]      \right\}\\
                     &=m^{3-\alpha}\left\{\sum_{k=4}^{\infty}{ 3-\alpha \choose k } \left[\left(\frac{2}{m}\right)^k  -4\left(\frac{1}{m}\right)^k -4\left(\frac{-1}{m}\right)^k +\left(\frac{-2}{m}\right)^k   \right]      \right\}\\
                     &=\frac{1}{m^{\alpha-1}}\left[  \frac{(3-\alpha)(2-\alpha)(1-\alpha)(-\alpha)}{m^2}+\cdots
                     \right]\\
                     &>0,
  \end{split}
  \end{equation*}
noting that the terms (with $k$ being odd numbers) of the series are
zero and the terms (with $k$ being even numbers) are positive.
According to Remark 2.1, all the corresponding properties of
$p_{i,k}^{\alpha}$ are obtained.
\end{proof}


\begin{thebibliography}{99}


\bibitem{Deng:07}  W.H. Deng, Numerical algorithm for the time fractional Fokker-Planck equation, J. Comput. Phys. 227 (2007) 1510-1522.

\bibitem{Diethelm:02} K. Diethelm, N.J. Ford, A.D. Freed, A predictor corrector
approach for the numerical solution of fractional differential
equations, Nonl. Dyna. 29 (2002) 3-22.

\bibitem{Gustafsson:95} B. Gustafsson, H.-O. Kreiss, J. Oliger, Time dependent problems and difference methods, John Wiley \& Sons, Inc, New York, 1995.

\bibitem{Kenneth:93}K.S. Miller, B. Ross, An Introduction to the Fractional Calculus and Fractional
 Differential Equations, Wiley-Interscience Publication, US, 1993.

\bibitem{Lapidus:82} L. Lapidus, G.F. Pinder, Numerical Solution of Partial Differential Equations in Science and Engineering, Wiley, New York, 1982.


\bibitem{Li:11} C.P. Li, A. Chen, J.J. Ye, Numerical approaches to fractional calculus and fractional ordinary
 differential equation, J. Comput. Phy. 230 (2011) 3352-3368.

\bibitem{Li:07}  C.P. Li, W.H. Deng, Remarks on fractional derivatives,
 Appl. Math. Comput. 187 (2007) 777-784.

 \bibitem{Liu:03} F. Liu, V. Anh, I. Turner, P. Zhuang, Time fractional advection
dispersion equation, J. Appl. Math. Computing 13 (2003) 233-245.

 \bibitem{Liu:07} F. Liu, P. Zhuang, V. Anh, I. Turner, K. Burrage, Stability and convergence of the
difference methods for the space-time fractional advection-diffusion equation, Appl.
Math. Comp., 191 (2007) 12-20.

\bibitem{Meerschaert:04} M.M. Meerschaert, C. Tadjeran, Finite difference approximations for fractional advection-dispersion
flow equations, J. Comput. Appl. Math. 172 (2004) 65-77.

\bibitem{Meerschaert:06} M.M. Meerschaert, C. Tadjeran, Finite difference approximations for two-sided
    space-fractional partial differential equations, Appl. Numer. Math. 56 (2006) 80-90.

\bibitem{Mark:06} M.M. Meerschaert, H.P.  Scheffler, C. Tadjeran, Finite difference methods for two-dimensional fractional dispersion equation,
      J. Comput. Phys. 211 (2006) 249-261.

\bibitem{Metzler:00} R. Metzler, J. Klafter, The random walk's guide to anomalous
 diffusion: A fractional dynamics approach, Phys. Rep. 339 (2000) 1-77.

\bibitem{Podlubny:99}I. Podlubny, Fractional Differential Equations, Academic Press, New York, 1999.

\bibitem{Samko:93} S. Samko, A. Kilbas, O. Marichev, Fractional Integrals and Derivatives: Theory and Applications, Gordon and Breach, London, 1993.

\bibitem{Sousa:Pro} E. Sousa, How to approximate the fractional derivative of order $1<\alpha \leq2 $, in: Proceedings
 of the 4th IFAC Workshop on Fractional Differentiation and its Applications, Badajoz, Spain, 2010.

\bibitem{Sousa:11} E. Sousa, Numerical approximations for fractional diffusion equations via splines,
 Comput. Math. Appl.  62  (2011)  938-944.

\bibitem{Sousa:12}   E. Sousa, C. Li, A weighted finite difference method for the fractional diffusion equation based on the
                    Riemann-Liouville drivative, arXiv:1109.2345v1 [math.NA].

\bibitem{Charles:07} C. Tadjeran, M.M. Meerschaert, A second-order accurate numerical method for the two-dimensional fractional diffusion equation,
       J. Comput. Phys. 220 (2007) 813-823.

\bibitem{Charles:06} C. Tadjeran, M.M. Meerschaert, H.P. Scheffler, A second-order accurate numerical approximation for
 the fractional diffusion equation, J. Comput. Phys. 213 (2006) 205-213.

\bibitem{Yang:10} Q. Yang, F. Liu, I. Turner, Numerical methods for fractional partial differential
equations with Riesz space fractional derivatives, Appl. Math. Model., 34 (2010) 200-218.


\end{thebibliography}
\end{document}